\documentclass[11pt]{article}
\usepackage{eurosym}
\usepackage{amsfonts}
\usepackage{amsmath}
\usepackage{geometry}
\usepackage{amssymb}
\usepackage{graphicx}
\usepackage{color}
\usepackage{tikz}
\usepackage[author={Oana}]{pdfcomment}
\usepackage[colorinlistoftodos]{todonotes}
\usepackage{titling}
\usepackage{lipsum}
\usepackage{mathrsfs} 
\usepackage{graphicx}
\usepackage{caption}
\usepackage{multicol}
\usepackage{lipsum}
\usepackage{mwe}
\usepackage{url,hyperref,lineno,microtype,subcaption}
\usepackage[onehalfspacing]{setspace}
\usepackage{subcaption,graphicx}
\usepackage{wrapfig}

\newtheorem{theorem}{Theorem}

\newtheorem{remark}[theorem]{Remark}

\newcommand{\bx}{\boldsymbol{x}}
\newcommand{\bmu}{\boldsymbol{\mu}}

\newcommand{\T}{\mathbb{T}}

\newcommand{\cL}{\mathcal{L}}
\newcommand{\cA}{\mathcal{A}}

\begin{document}

\date{}
\title{\textbf{Bayesian Inference for Fluid Dynamics: A Case Study for the Stochastic Rotating Shallow Water Model}}
\author{ Peter Jan van Leeuwen \qquad Dan Crisan \qquad Oana Lang \qquad Roland Potthast \\
}
\maketitle

\begin{abstract}
    In this work, we use a tempering-based adaptive particle filter to infer from a partially observed stochastic rotating shallow water (SRSW) model which has been derived using the Stochastic Advection by Lie Transport (SALT) approach introduced in \cite{Holm2015}.  The methodology we present here validates the applicability of tempering and sample regeneration using a  Metropolis-Hastings methodology  to high-dimensional models used in stochastic fluid dynamics. The methodology is first tested on the Lorenz 63 model with both full and partial observation. We study the efficiency of the particle filter for the Lorenz 63 model as well as the  SRSW model. 
    
\end{abstract}

\section{Introduction}\label{particlefiltersection}

Let $X$ and $Z$ be two processes defined on the probability space $(\Omega, \mathcal{F}, \mathbb{P})$. $X$ is usually called the \textit{signal process} or the \textit{truth} and $Z$ is the \textit{observation process}. In our case $X$ is given by a pathwise solution of the stochastic rotating shallow water system computed on a staggered grid. We have proven in \cite{CL3} that such a solution exists.\footnote{Note that in \cite{CL3} we work with infinite dimensional function spaces, while here the state space will be $\mathbb{R}^{d_X}$.} The nonlinear filtering problem consists in finding the best approximation of the posterior distribution of the signal $X_t$ given the observations $Z_1, Z_2, \ldots, Z_t$ \footnote{For a mathematical introduction on the subject see \cite{CrisanBain}. For an introduction from the data assimilation perspective see \cite{Leeuwenbook2015} and \cite{ReichCotter2015}.}. The posterior distribution is usually denoted by $\pi_t$. We denote the dimension of the state space by $d_X$ and the dimension of the observation space by $d_Z$.  

\noindent A \textit{particle filter} is a sequential Monte Carlo method in which the posterior is approximated using a set of \textit{particles}, that is random measures of the form
\begin{equation*}
\pi_t \approx \displaystyle\sum_{\ell} \mathrm{w}_{t}^{\ell}\delta({x_{t}^{\ell}})
\end{equation*}
where $\delta$ is the Dirac delta function, $\mathrm{w}_{t}^1, \mathrm{w}_{t}^2, \ldots $ are the \textit{weights} of the  particles and $x_{t}^1, x_{t}^2, \ldots$ are their corresponding positions \cite{CrisanBain}.
 One can make inferences about the signal process using Bayes' theorem, the time-evolution induced by the model, and observations \cite{ReichCotter2015}, \cite{CrisanBain}. Observations are modelled as noisy measurements of the truth, using the \textit{observation operator}:
 \begin{equation}\mathscr{H}: \mathbb{R}^{d_X} \rightarrow \mathbb{R}^{d_Z}\end{equation}
  \begin{equation}Z_t = \mathscr{H}(X_t)+ V_t\end{equation}
  where $(V_t)_{t\geq 0}$ are independent identically distributed random variables with standard normal distribution and $\mathscr{H}$ is a Borel-measurable function. 
 Observations are incorporated into the system at \textit{assimilation times}. The ensemble of particles is evolved between assimilation times according to the law of the signal. At each assimilation time the observation is incorporated into the system through the \textit{likelihood function}: 
 \begin{equation}
 g_t^{z_t}:\mathbb{R}^{d_X} \rightarrow [0,1], \ g_t^{z_t}(x) = g_t(z_t - \mathscr{H}(x_t)) = \mathbb{P}(Z_t \in dz_t | X_t = x_t)\end{equation} 
 that is 
 \begin{equation}
 \displaystyle\int_A g(z_t-\mathscr{H}(x_t))dz_t = \mathbb{P}(Z_t \in A|X_t=x_t)
 \end{equation}
 where $A \in \mathcal{B}(\mathbb{R}^{d_Z})$. 
 The following recursion formula holds (see \cite{CrisanBain})
 \begin{equation}\label{recursionformula}
 \pi_t = g_t \star \pi_{t-1}\mathcal{K}_t
 \end{equation}
 where  
 \begin{equation} \mathcal{K}_t:\mathbb{R}^{d_X}\times\mathcal{B}(\mathbb{R}^{d_X}) \rightarrow [0,1], \ \mathcal{K}_t(X_{t-1}, B) = \mathbb{P}(X_t \in B|X_{t-1})
 \end{equation} 
 for any measurable set $B \in \mathcal{B}(\mathbb{R}^{d_X})$ and by '$\star$' we denoted the projective product as defined in the Appendix. 
 Schematically, this can be arranged as  
 \begin{equation}
\pi_{t-1}^{a, z_{0:t-1}} \xrightarrow[\substack {model\\forecast\\prediction}]{\mathcal{K}_t} \pi_{t-1}^{a, z_{0:t-1}}\mathcal{K}_t =: \pi_{t}^b =: p_t \xrightarrow[\substack {assimilation\\analysis \\update}]{{tempering}, \ g_t^{z_t} \star} g_t^{z_t} \star\pi_{t}^b=\pi_t^{a, z_{0:t}}. 
 \end{equation} 
 The indices $a$ and $b$ stand for \textit{analysis} and \textit{background} respectively. Here $p_t=p_t^{Z_{0:t-1}}$ is the prior distribution of the signal, that is a probability measure which belongs to the space $ \mathcal{P}(\underbrace{\mathbb{R}^{d_X}\times\mathbb{R}^{d_X}\times \ldots \times \mathbb{R}^{d_X}}_{t+1})$ with $t=0,1, \ldots$, and corresponds to the path space generated by $X_0, X_1, \ldots X_t$. By $Z_{0:t}$ we mean the random vector $(Z_0, Z_1, \ldots, Z_{t})$ with the corresponding observed values $z_{0:t}$. Taking into account the definition of the projective product, the recursion formula \eqref{recursionformula} can also be written as
 \begin{equation}
 \pi_t(B) = \frac{\displaystyle\int_{B}g_t^{z_t}(x_t)p_t(dx_t)}{\displaystyle\int_{\mathbb{R}^{d_X}}g_t^{z_t}(x_t)p_t(dx_t)} = {\beta}_t^{-1}\displaystyle\int_{B}g_t^{z_t}(x_t)p_t(dx_t)
 \end{equation}
 where $B \in \mathcal{B}(\mathbb{R}^{d_X})$ and $\beta_t:= \displaystyle\int_{\mathbb{R}^{d_X}}g_t^{z_t}(x_t)p_t(dx_t)$ is a normalising constant. 
  We are looking for an approximation of the form 
   \begin{equation}
   \pi_t \approx \pi_t^N = \displaystyle\sum_{\ell=1}^N \mathrm{w}_t^{\ell}\delta(x_t^{\ell}).
    \end{equation}
   The typical particle filter runs as follows. 
   Each particle is propagated forward in time using the model (SRSW or Lorenz '63 in our case). We model the evolution of the signal discretely in time through a map $\mathcal{M}_t$. In the following, $\mathcal{M}_t$ is a discrete realisation of the Lorenz '63 or the SRSW model. We have
    \begin{equation} 
    \mathcal{M}_{t}: \mathbb{R}^{d_X} \rightarrow \mathbb{R}^{d_X}
    \end{equation}
   and if $x_{t_1}^{\ell}$, $x_{t_2}^{\ell}$ is the position of the particle $\ell$ at time $t_1$ respectively $t_2$ then
    \begin{equation} 
    x_{t_2}^{\ell} = \mathcal{M}_{t_2}(x_{t_1}^{\ell}).
    \end{equation} 
   The particle trajectory is independent of the trajectory of the signal. 
   At each assimilation time, every particle is weighted depending on the likelihood of its position, given the observation. The weight measures how close the particle trajectory is to the signal trajectory. 
   \begin{equation} 
   \pi_{t_2}^a = \displaystyle\sum_{\ell=1}^N \mathrm{w}_{t_2}^{\ell}\delta(\mathcal{M}_{t_2}(x_{t_1}^{\ell})) = \displaystyle\sum_{\ell=1}^N \mathrm{w}_{t_2}^{\ell} \delta(x_{t_2}^{\ell})
   \end{equation} 
  Particles which are close to the \textit{truth} and therefore have higher weights will be multiplied, while those which are far away will be eliminated. There are several studies (\cite{VetraCarvalho-Leeuwen2018},\cite{Leeuwen.et.al.2019}) which have shown that in high-dimensional spaces the tendency is to have one particle gaining a weight close to one while all the others are discarded, having weight close to zero. Indeed, in \cite{CrisanBain} it has been proven that a particle filter will provide a good approximation of the posterior distribution only when enough particles are used, which explains the difficulty of the problem at hand. The standard rapid divergence of the particles' trajectory from the signal trajectory is known in the literature as \textit{the curse of dimensionality} \cite{VetraCarvalho-Leeuwen2018}, \cite{Leeuwenbook2015}. Nonetheless, considerable progress has been made in this direction over the last years. A state-of-the-art analysis of the most recent efforts on tackling the filter degeneracy problem can be found in \cite{Leeuwen.et.al.2019} and \cite{VetraCarvalho-Leeuwen2018}. Innovative remedies arise from different directions: optimal transportation \cite{ReichCotter2015}, tempering \cite{CrisanBeskosJasra},\cite{KantasBeskosJasra},\cite{Wei1},\cite{Wei2}, localisation \cite{ReichChen},\cite{Potthast}, model reduction \cite{Wei1}, data assimilation as a boundary value problem \cite{ReichActaNumerica}, jittering \cite{CrisanBeskosJasra}, nudging \cite{Shevchenko.et.al}, and proposal densities \cite{Leeuwen2010}. Some of them have been tested in operational numerical weather prediction systems, e.g. \cite{Potthast}. 
 
  Our work is part of these attempts of developing a particle filter methodology for high dimensional models originating in (stochastic) fluid dynamics. We have implemented a particle filter with adaptive tempering and jittering for the stochastic rotating shallow water model. The suitability of these two methods used in tandem for high-dimensional problems has been proven in \cite{CrisanBeskosJasra} and tested in \cite{Wei1}, \cite{Wei2}, \cite{Shevchenko.et.al}. In \cite{Wei2} the method is used for the stochastic incompressible 2D Euler model with damping and forcing, while in \cite{Shevchenko.et.al} it is tested for the 2D quasi-geostrophic model. We have no knowledge of any previous application of this methodology for the stochastic rotating shallow water system. The complexity of this model makes the problem harder. The aim of this paper is to show the applicability of the particle filter described below and the stochastic rotating shallow water model running alongside. We start with first coupling the particle filter with the classical Lorenz '63 model and then we move to the targeted SRSW model. The results are positive in both cases.  
 
 \section{The Algorithm}
 
  In tempering, an artificial dynamics is introduced via a sequence of artificial target distributions between prior and posterior. Each intermediate distribution has a characteristic \textit{temperature} chosen such that a reasonable number of good particles survive. We first follow a resampling procedure \footnote{This is a static procedure.} in which the particles with low weights are replaced with particles with higher weights such that at the end an ensemble of equal-weighted particles is obtained. However, the particles can be concentrated in the wrong place. In order to quantify the spread of the weights with respect to the posterior, we use the \textit{effective sample size} statistic:
  
  \begin{equation} 
  ess(\mathrm{w}) = \frac{1}{\displaystyle\sum_{\ell=1}^N (\mathrm{w}^{\ell})^2}
  \end{equation}
  
  The \textit{ess} will be chosen to be above a certain threshold in order for the ensemble of particles to be oriented in the right direction. In tempering one increases gradually the variance of the distribution 
  using a sequence of \textit{temperatures} $0=\phi_0 < \phi_1 < \ldots \phi_R=1$ to ensure that the $ess$ remains above the chosen threshold. Once the temperature is (dynamically) chosen, a resampling procedure is applied. The output is a sequence of \textit{tempered posterior distributions} with corresponding normalised tempered weights. In this case the intermediate tempered posterior distribution at the $r^{th}$ tempering step is given by 
     \begin{equation}
     \pi_t^r(B) = \frac{\displaystyle\int_{B}\left(g_t^{z_t}(x_t)\right)^{\phi_r}p_t(dx_t)}{\displaystyle\int_{\mathbb{R}^{d_X}}\left(g_t^{z_t}(x_t)\right)^{\phi_r}p_t(dx_t)}
     \end{equation}
  for any $B \in \mathcal{B}(\mathbb{R}^{d_X})$. We denote by $\textbf{x}=(x^{\ell})_{\ell=1}^N$ the ensemble of particles. Then
  \begin{equation}
  \pi_{t_i}^{r,N} = \displaystyle\sum_{\ell=1}^{N} \mathrm{w}_{t_i}^{r,\ell}(\phi_r,\textbf{x})\delta(x_{t_i}^{\ell})
  \end{equation}
  where 
          \begin{equation} 
          \mathrm{w}_{t_i}^{r,\ell}(\phi_r, \textbf{x}) = \left(g_{t_i}^{z_{t_i}}(x_{t_i}^{\ell})\right)^{\phi_r-\phi_{r-1}}, \ \ \ \hbox{with} \ \ \ \displaystyle\sum_l \mathrm{w}_{t_i}^{\ell} = 1.
          \end{equation}
     The corresponding $ess$ is given by
    \begin{equation}
    ess_i(\phi^r, \textbf{x}) := \|\mathrm{w}_{t_i}(\phi_r, \textbf{x})\|_{{\ell}^2}^{-1}.
    \end{equation} 

  The following is the tempering and jittering algorithm (see e.g. \cite{KantasBeskosJasra},\cite{Wei1}):
  \begin{itemize}
  {\item [1.] At initial time $t=0$: sample $N$ particles from the prior distribution.}
  \vspace{1mm} 
  { \item[2.] On the time interval $(t_{i-1}, t_i]:$ we have an ensemble $\textbf{x}$ of particles with positions $(x_{t_{i-1}}^{\ell})_{\ell}$ and we want to assimilate observational data $z_{t_i}$ in order to obtain a new ensemble $(x_{t_i}^{\ell})_{\ell}$ that defines $\pi_{t_i}^N$:}
     \begin{itemize}
         \item[2.1.] Evolve $x_{t_{i-1}}^{\ell} \xrightarrow[\substack {SRSW, Lorenz63}]{SPDE} x_{t_i}^{\ell}$.
         \vspace{1mm} 
        \item[2.2.] Set temperature $\phi = 1$.
        \vspace{1mm} 
        \item[2.3.] While $ess_i(\phi, \textbf{x}) < N_{threshold}$ do
            \begin{itemize}
                \item Find $\phi' \in (1-\phi,1)$ such that $ess_i(\phi'-(1-\phi), \textbf{x}) \approx N_{threshold}$. Resample according to $\mathrm{w}_{t_i}^{\ell}(\phi'-(1-\phi), \textbf{x})$ and apply MCMC with jittering if required (i.e. if there are duplicates) $\Rightarrow$ a new ensemble $\textbf{x}(\phi')$.
               \vspace{1mm} 
               \item Set $\phi = 1-\phi'$ and $\textbf{x}=\textbf{x}(\phi)$.
            \end{itemize}
        \item[2.4.] If $ess_i \geq N_{threshold}$ then Stop and go to the $(i+1)^{th}$ filtering step with $(x_{t_i}^{\ell}, \mathrm{w}_i^{\ell})_{\ell}$.
     \end{itemize}  
  \end{itemize}
The jittering procedure is standard (see e.g. \cite{Wei1}) but we will briefly explain the idea. Without jittering we have 
\begin{equation}
x_{t_2}^{\ell} = \mathcal{M}(x_{t_1}^{\ell})
\end{equation}
which can also be written as
\begin{equation}
x_{t_2}^{\ell} = \mathcal{M}(x_{t_1}^{\ell}, W(t_1:t_2))
\end{equation}
where $W$ is the driving Brownian motion of the model and with $W(t_1:t_2)$ we denoted the Brownian path between $t_1$ and $t_2$. In this case the prior distribution is given by
\begin{equation}
p_{t_2} = \frac{1}{N}\displaystyle\sum_{\ell = 1}^{N} \delta (x_{t_2}^{\ell}).
\end{equation} 
The problem is that, when the number of independent observations is large, after resampling, the particles end up in the same place and we can have a large number of duplicates. In order to overcome this issue we modify the last part of the particle trajectory using a \textit{jittering parameter} $\rho$ and a Gaussian random variable $Z$ which is orthogonal to $W$. Then the new dynamics is given by
\begin{equation}\label{rhoeqn}
\tilde{x}_{t_2}^{\ell} = \mathcal{M}_{t_2}\left(x_{t_1}^{\ell}, \rho W(t_1:t_2) + \sqrt{1-\rho^2}Z(t_1:t_2)\right)
\end{equation}
where
\begin{equation}
x_{t_i}^{\ell} = \mathcal{M}_{t_i}\left(x_{t_{i-1}}^{\ell}, W(t_{i-1}:t_i)\right).
\end{equation}

\section{Applications for the Lorenz 1963 Model}
\subsection{Model Description}
The \textit{Lorenz '63} model is a classical nonlinear three-dimensional model that is a precursor of turbulence theory which has been introduced in \cite{Lorenz63}. It reads
\begin{subequations}\label{lorenz}
\begin{alignat}{3} 
&\frac{dx}{dt} = \alpha (y-x) \label{lorenz1}\\
&\frac{dy}{dt} = (\beta - z)x - y \label{lorenz2}\\
&\frac{dz}{dt} = xy - \gamma z \label{lorenz3}
\end{alignat} 
\end{subequations}
where $\alpha, \beta, \gamma$ are real positive parameters. The model is well-known for the broad spectrum of patterns displayed for different values of $\alpha, \beta, \gamma$ and its well-known \textit{butterfly attractor}. The original values chosen by Lorenz in \cite{Lorenz63} were $\alpha=10, \beta = 28$ and $\gamma = \frac{8}{3}$. For a discussion on the behaviour of the solutions for different parameter values see e.g. \cite{Sparrow}. 
This model is implemented using a Runge-Kutta scheme of order 4, with initial conditions $x_0 =1.508870, y_0=-1.531271, z_0=25.46091$. The details of the Runge-Kutta scheme can be found in Appendix. { The Lorenz '63 model has noise in it: after using the Runge-Kutta scheme to implement the three variables from system \eqref{lorenz}, we generate a random field and perturb the system in a manner which is similar to the one explained in \eqref{rhoeqn}. The parameter $\rho$ is equal to $0.99$ here.} 

\subsection{Data Assimilation Results}
\begin{figure}[ht!]
  \centering
  \hspace{-4em}
  \begin{subfigure}{.5\linewidth}
    \centering
    \includegraphics[width = \linewidth]{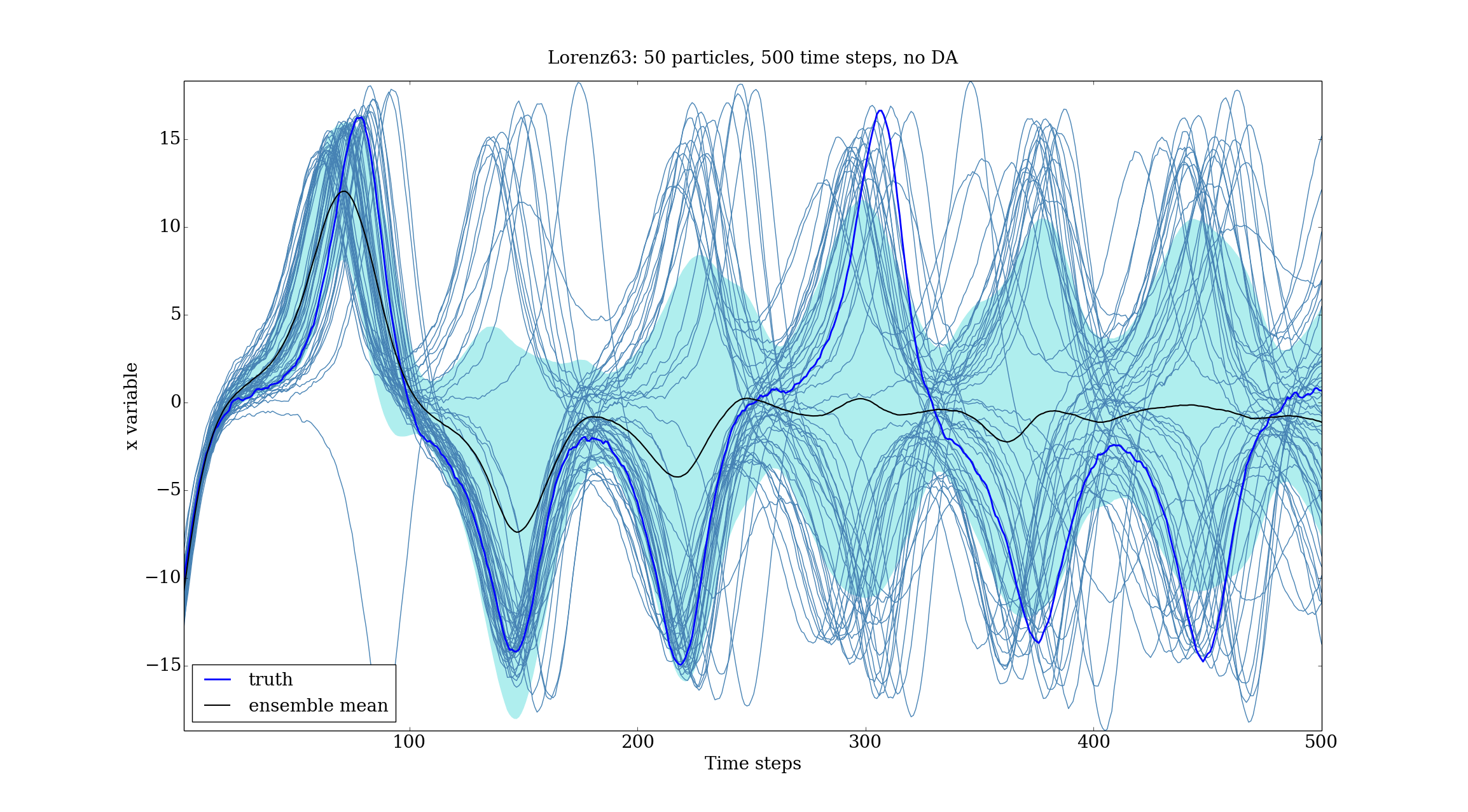}
    \caption{$x$ variable}
    \label{fignoDAx}
  \end{subfigure}%
  \begin{subfigure}{.5\linewidth}
    \centering
   \includegraphics[width = \linewidth]{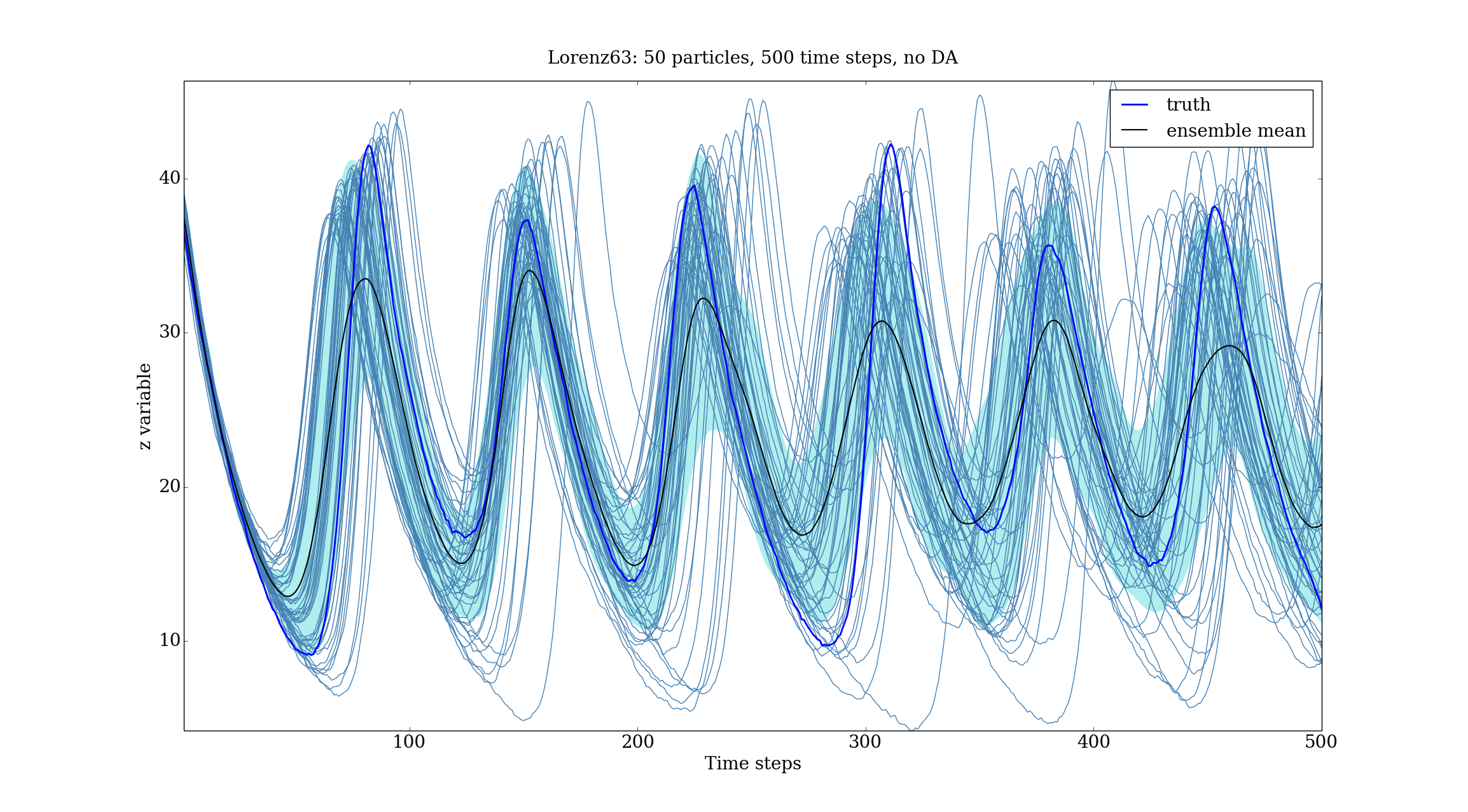}
   \caption{$z$ variable}
   \label{fignoDAz}
  \end{subfigure}%

  \caption{Evolution of the Lorenz '63 model for 500 time steps without any data assimilation.}
\end{figure}

\begin{figure}[ht!]
  \centering
  \begin{subfigure}{.5\linewidth}
    \centering
    \includegraphics[width = \linewidth]{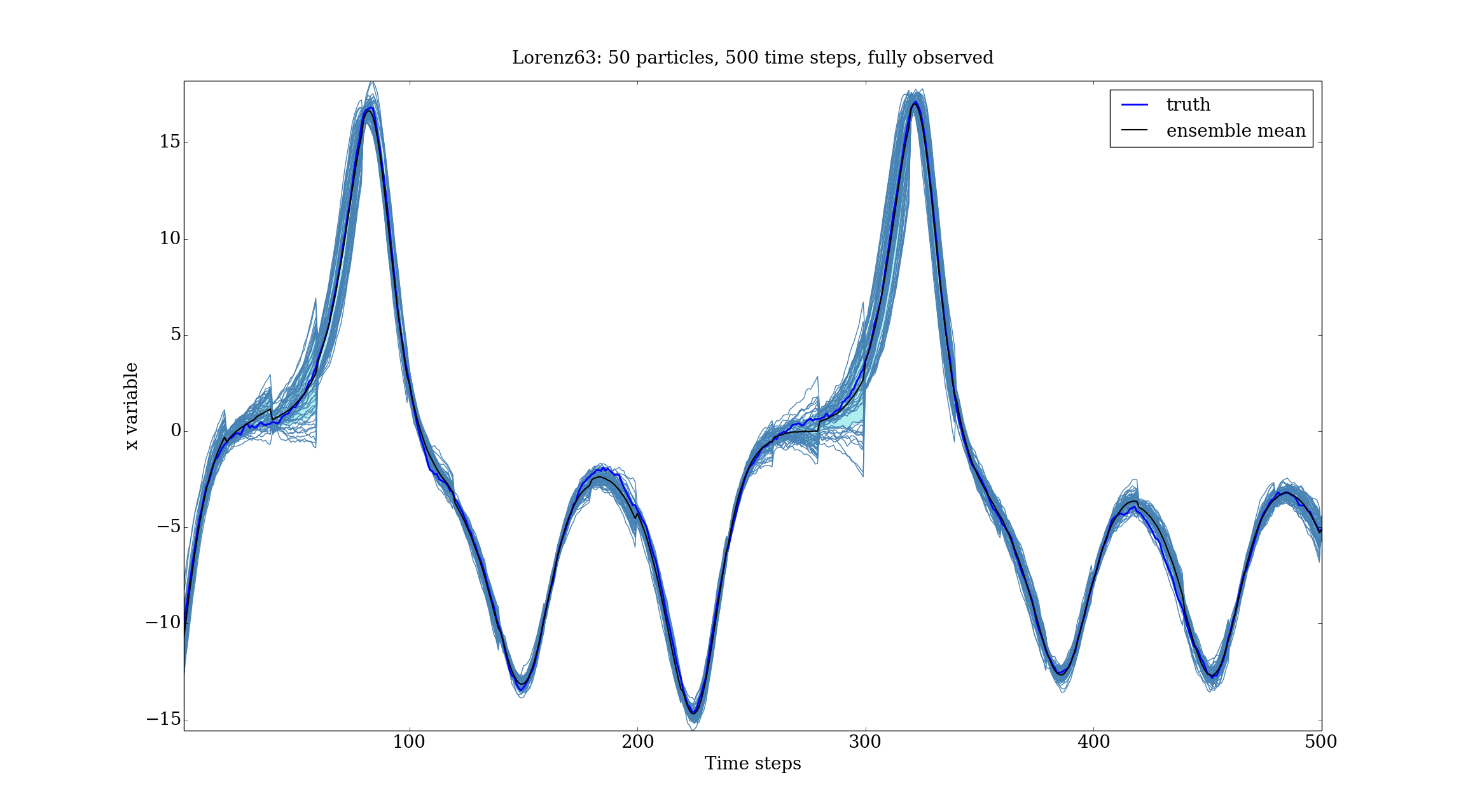}
    \caption{$x$ variable}
    \label{daS1x}
  \end{subfigure}%
  \begin{subfigure}{.5\linewidth}
    \centering
   \includegraphics[width = \linewidth]{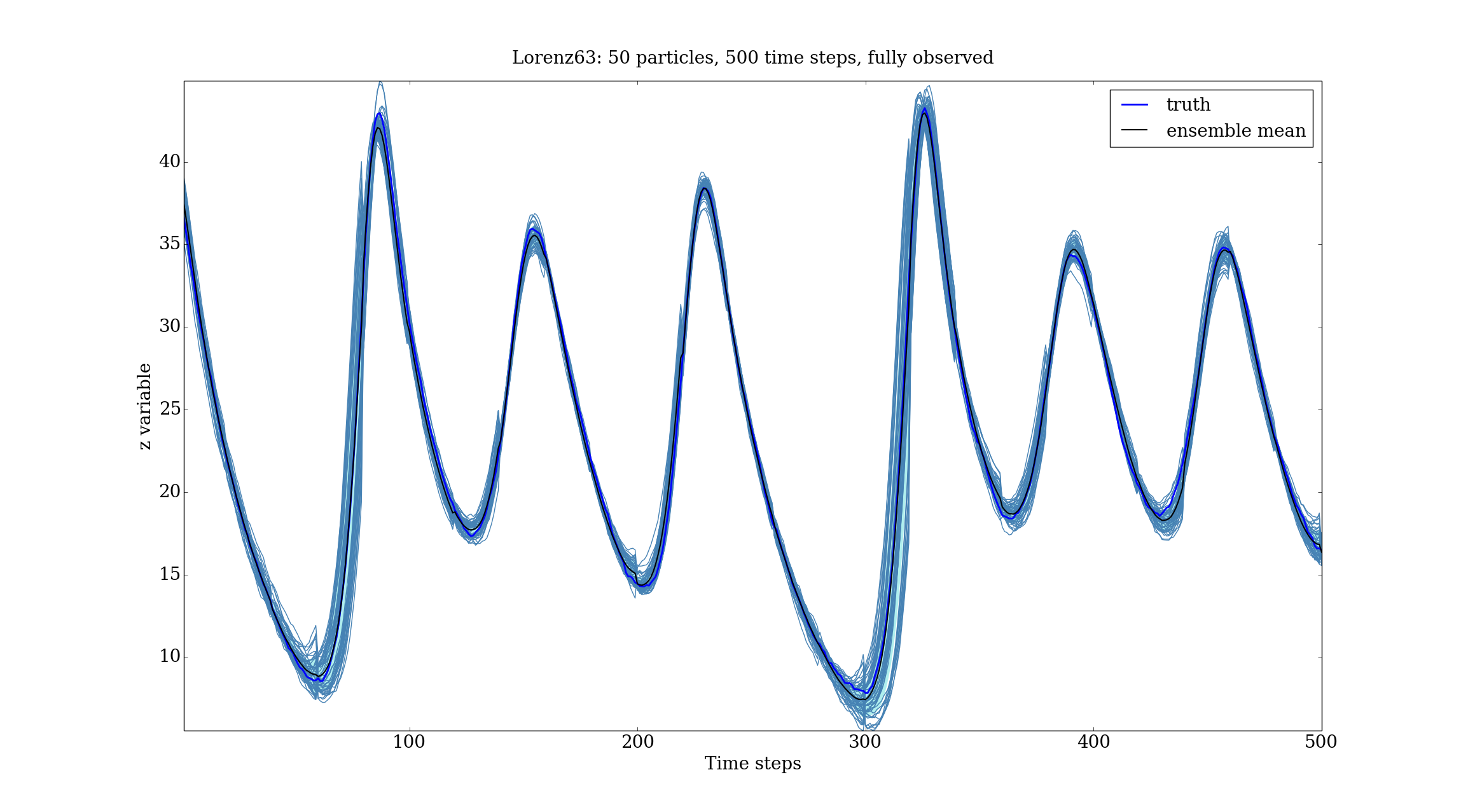}
   \caption{$z$ variable}
   \label{daS1z}
  \end{subfigure}%
  \caption{Evolution of the Lorenz '63 model for 500 time steps, all 3 variables are observed every 20 time steps. The initial uncertainty and the observational uncertainty are both equal to 1, the model error is equal to 0.1. Uncertainty is substantially reduced at assimilation times.}
\end{figure}

\begin{figure}[ht!]
  \centering
  \begin{subfigure}{.5\linewidth}
    \centering
    \includegraphics[width = \linewidth]{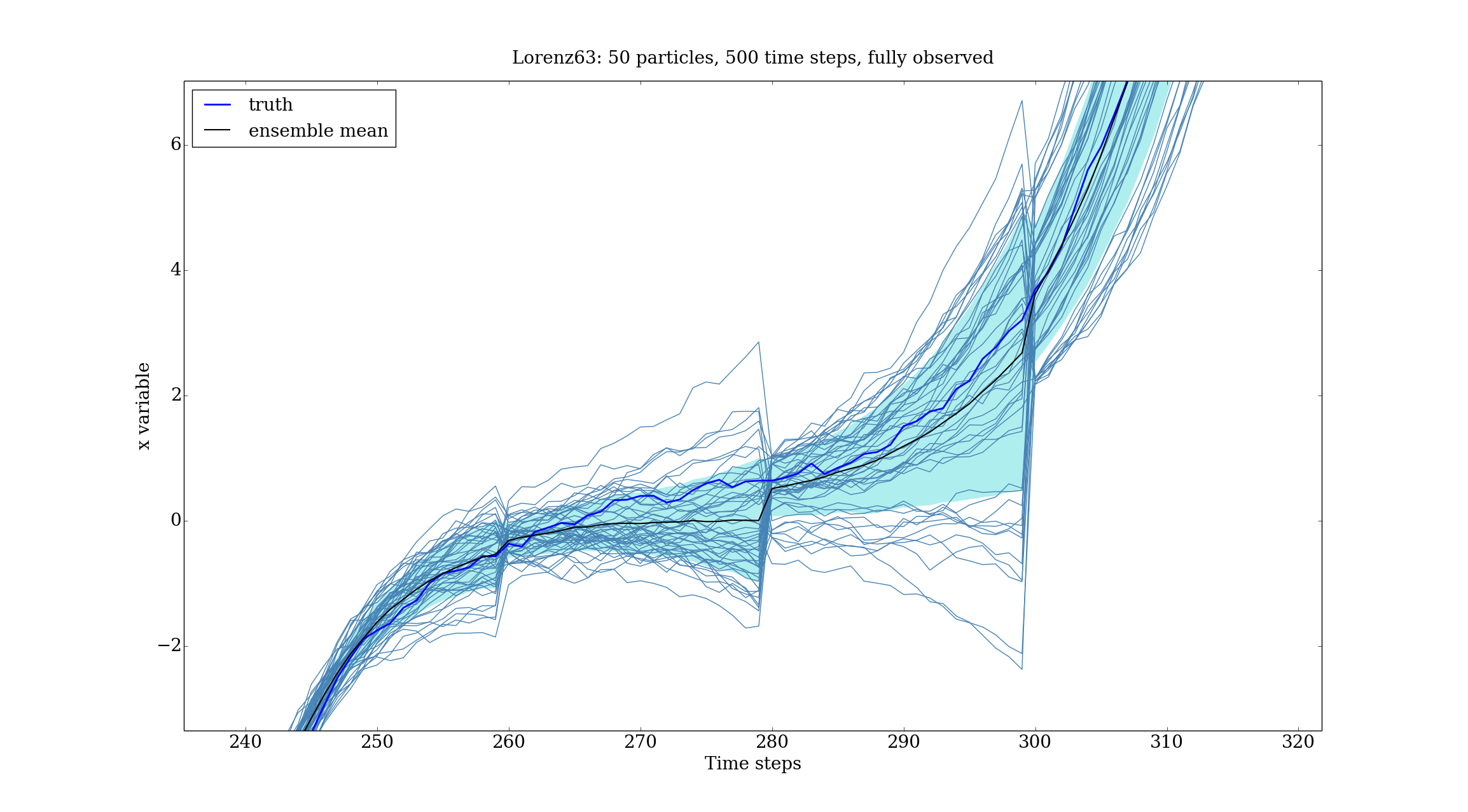}
    \caption{$x$ variable enlarged}
   \label{daS1enlarged}
  \end{subfigure}%
  \begin{subfigure}{.5\linewidth}
    \centering
   \includegraphics[width = \linewidth]{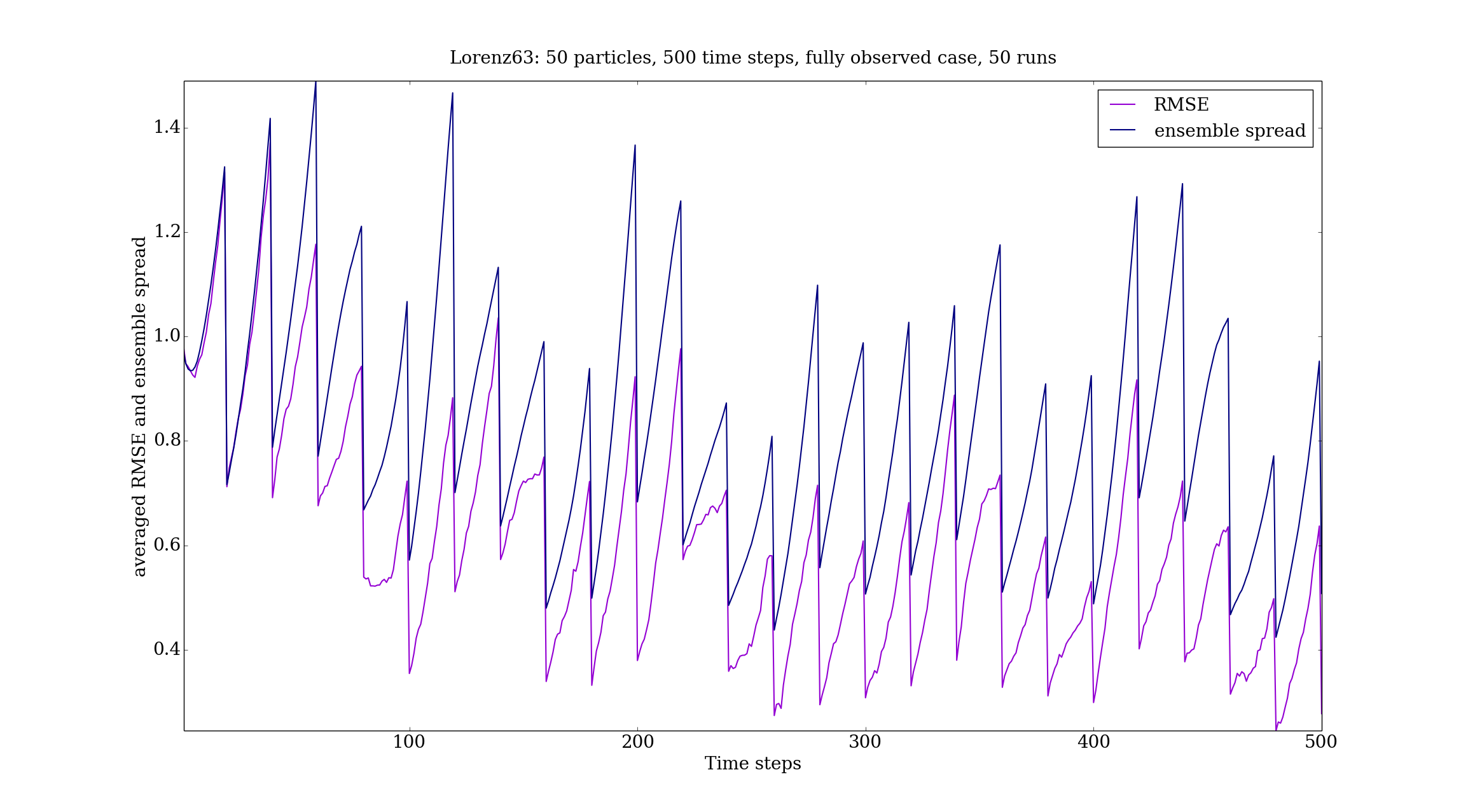}
   \caption{RMSE and ES}
   \label{daS1rmse}
  \end{subfigure}%
  \caption{Evolution of the Lorenz '63 model for 500 time steps, all 3 variables are observed every 20 time steps. The initial uncertainty and the observational uncertainty are both equal to 1, the model error is equal to 0.1. Uncertainty is substantially reduced at assimilation times.}
\end{figure}

We perform the data assimilation analysis using an ensemble of 50 particles which evolve for 500 time steps and each time step has size 0.01. In the standard scenario {(Figures \ref{daS1x}-\ref{daS1rmse})} all three variables of the system are observed every 20 time steps. The initial uncertainty is equal to 1, while the observational uncertainty and the model error are equal to 0.1. We present below a couple of scenarios obtained for different values of these parameters. The output is displayed for the first and the third variable. The ensemble of particles is plotted one standard deviation region around the ensemble mean. 
We first show {in Figures \ref{fignoDAx}-\ref{fignoDAz}} how the model evolves without any assimilation of data . 

\begin{wrapfigure}{l}{0.45\textwidth}
    \centering
    \includegraphics[width = \linewidth]{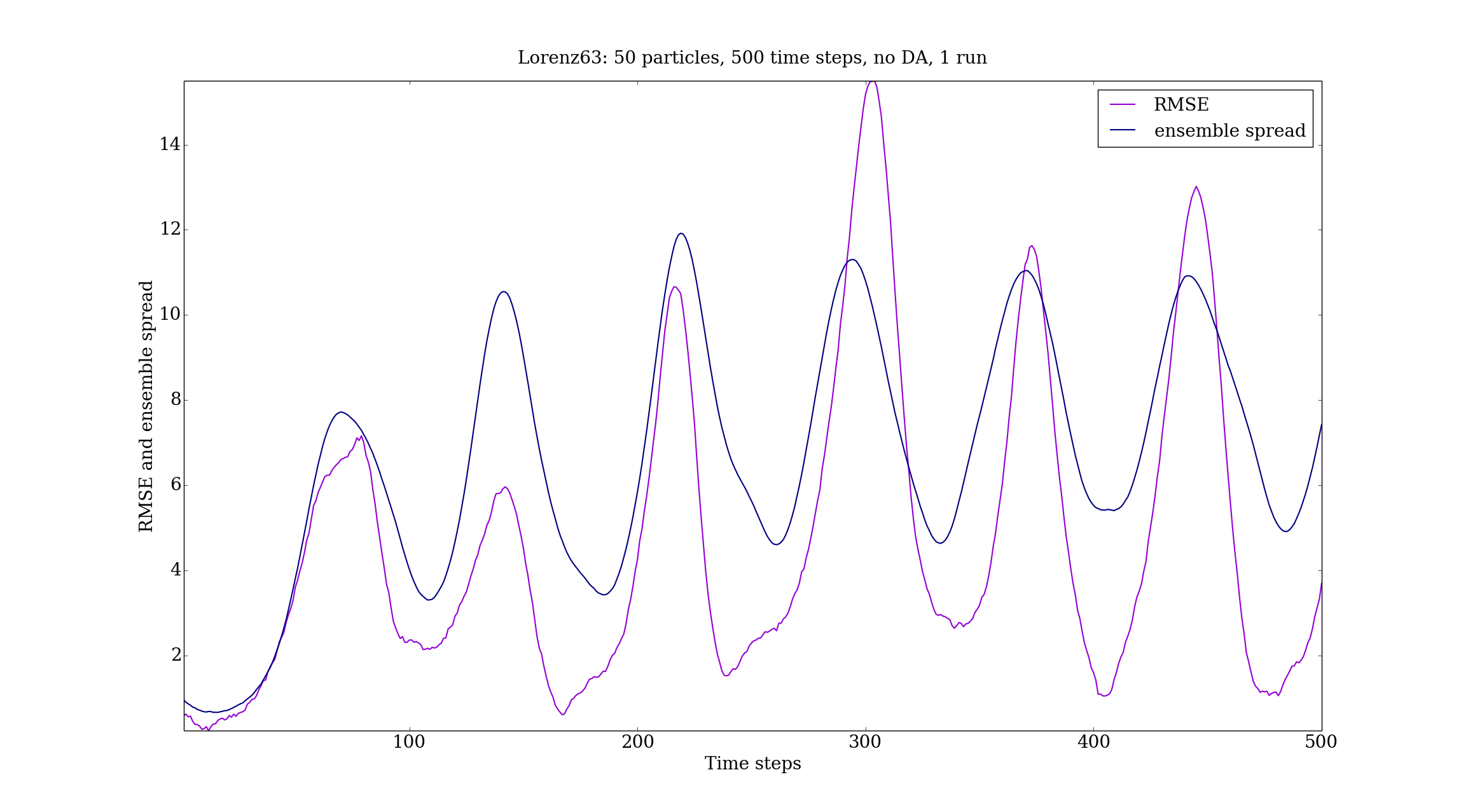}
    \caption{RMSE and ES}
    \label{fignoDArmse}
  \caption{Evolution of the Lorenz '63 model for 500 time steps without any data assimilation: RMSE and ES.}
\end{wrapfigure}

As it can be seen in Figure \ref{fignoDAx} and \ref{fignoDAz}, in the absence of the data assimilation step, the particles spread gradually around the entire attractor. This behaviour is exhibited when plotting both the first variable $x$ and the third variable $z$. In other words, the uncertainty increases (at times quite dramatically). The root mean square error (RMSE) and the ensemble spread (ES) is plotted in Figure \ref{fignoDArmse}. In this case the RMSE and the ES oscillate around a value which is comparable with the size of the attractor (the particles \textit{fill out} the attractor). 

Compared to the previous scenario, we show { in Figures \ref{daS1x}-\ref{daS1rmse}} that by assimilating observations every 20 time steps, we manage to substantially reduce the uncertainty. 
Through a repeated application of the forecast/data assimilation step, the cloud of particles successfully tracks the truth, evolving around the attractor set. 


In Figure \ref{daS1enlarged} we exhibit the first three instances to emphasize the reduction of uncertainty resulting from the application of the particle filter with tempering and jittering, as described in Section \ref{particlefiltersection}. 

We plot in Figure \ref{daS1rmse} the RMSE and the ES. We can see that the RMSE and the ES are comparable. This is a feature highly appreciated by data assimilation practitioners because it shows that our estimate of the uncertainty as measured by the width of the ensemble is a good estimate of the actual error in the ensemble mean. Based on this measure of success we can conclude that the particle filter performs satisfactory in the case where all the three variables of the system are observed. 

\begin{figure}[ht!]
  \centering
  \begin{subfigure}{.5\linewidth}
    \centering
    \includegraphics[width = \linewidth]{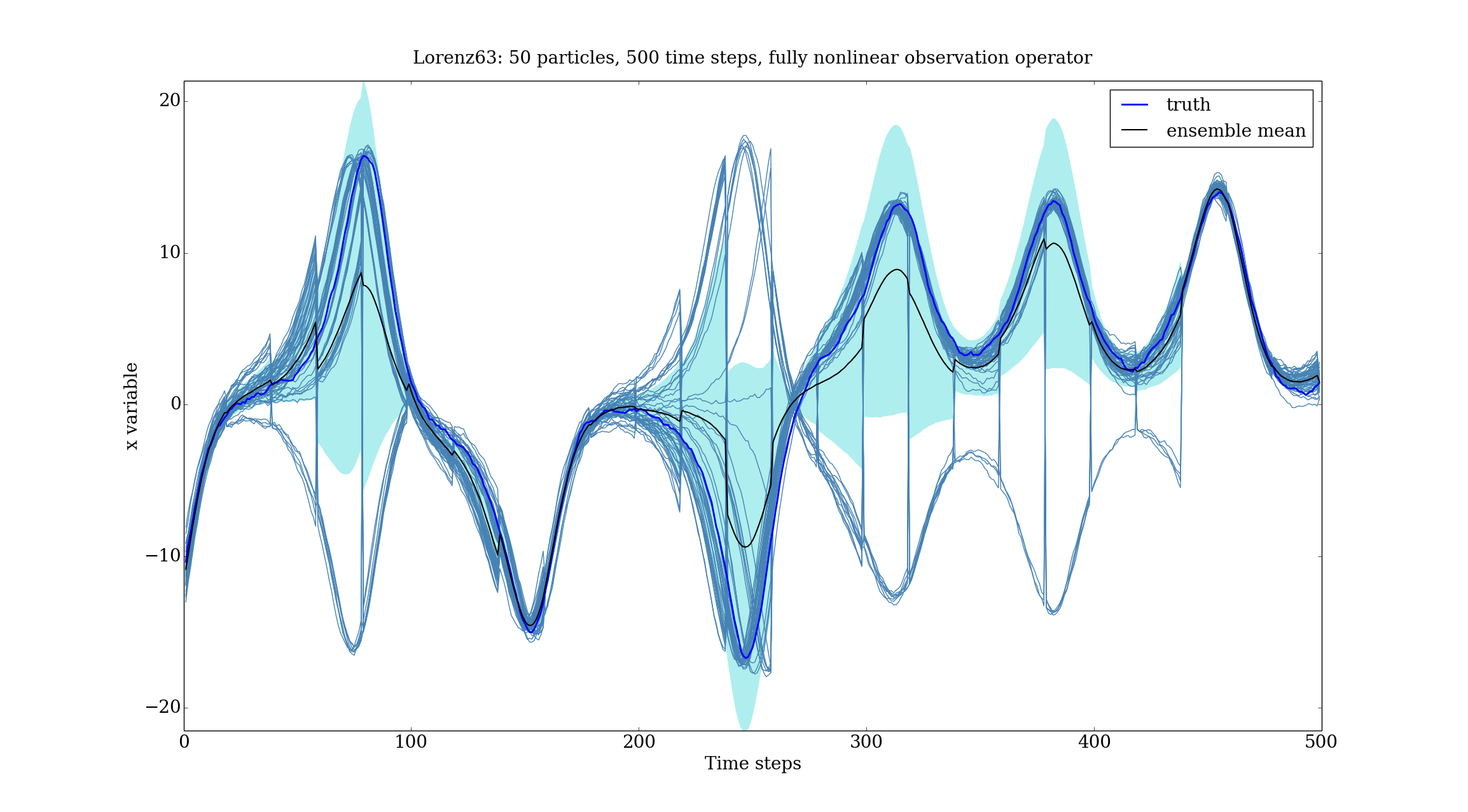}
    \caption{$x$ variable}
    \label{daS9x}
  \end{subfigure}%
  \begin{subfigure}{.5\linewidth}
    \centering
   \includegraphics[width = \linewidth]{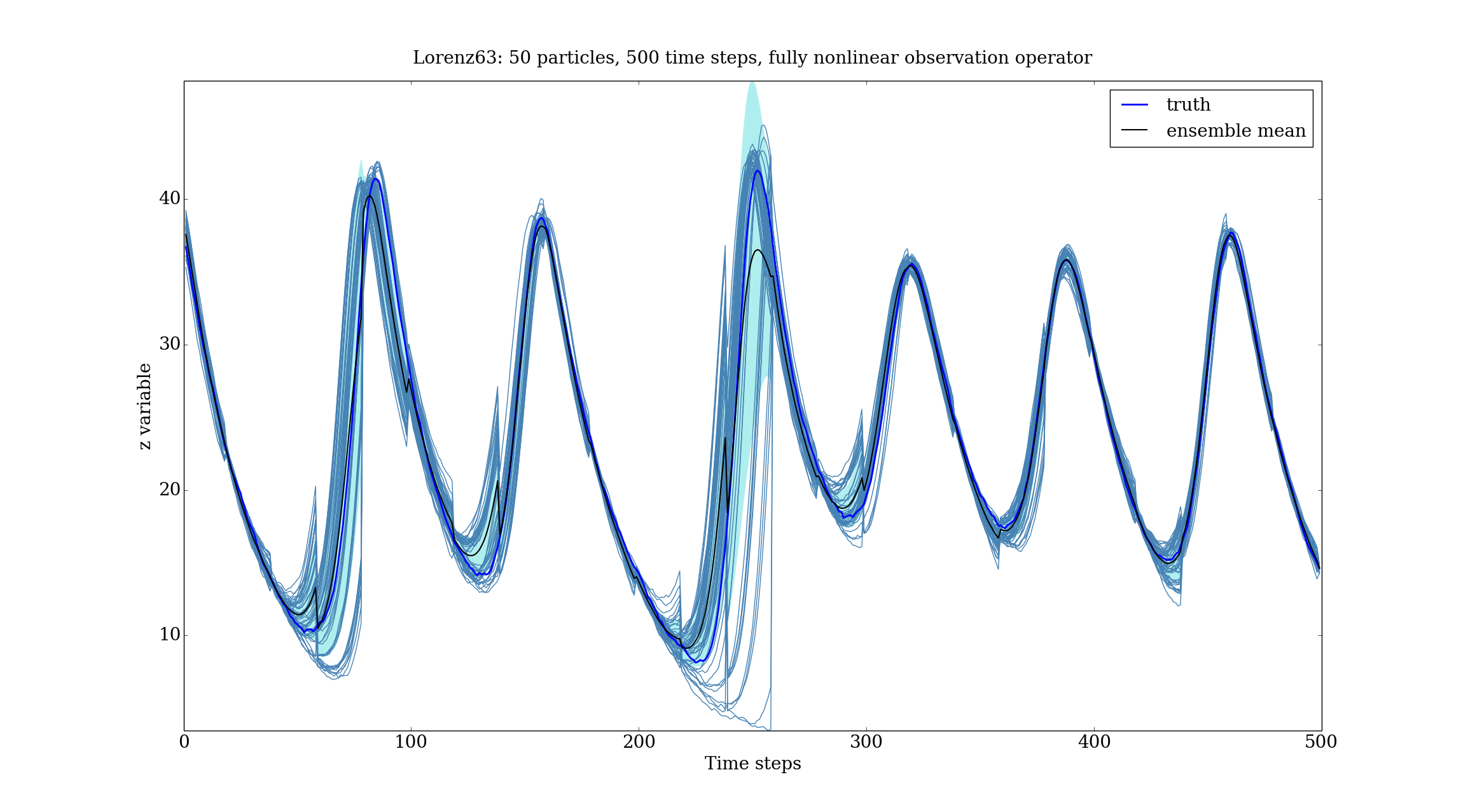}
   \caption{$z$ variable}
   \label{daS9z}
  \end{subfigure}%
  \caption{Evolution of the Lorenz '63 model for 500 time steps, all 3 variables are observed every 20 time steps. Observations are denoted by $z_1,z_2,z_3$. Fully nonlinear observation operator i.e. $z_1=x^2, z_2=y^2, z_3=z^2$ plus noise.}
\end{figure}

\begin{figure}[ht!]
  \centering
  \begin{subfigure}{.5\linewidth}
    \centering
    \includegraphics[width = \linewidth]{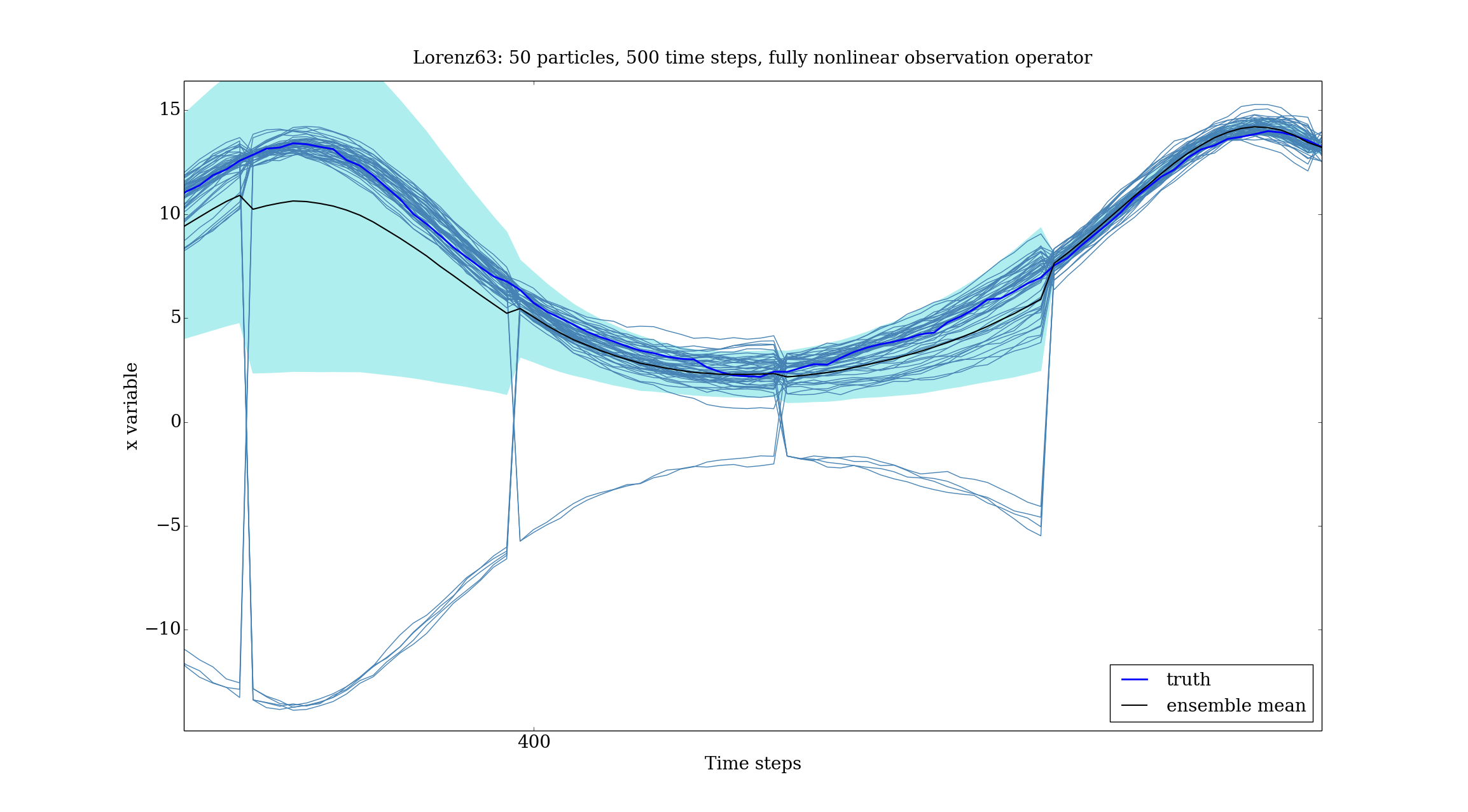}
    \caption{$x$ variable enlarged}
   \label{daS9enlarged}
  \end{subfigure}%
  \begin{subfigure}{.5\linewidth}
    \centering
   \includegraphics[width = \linewidth]{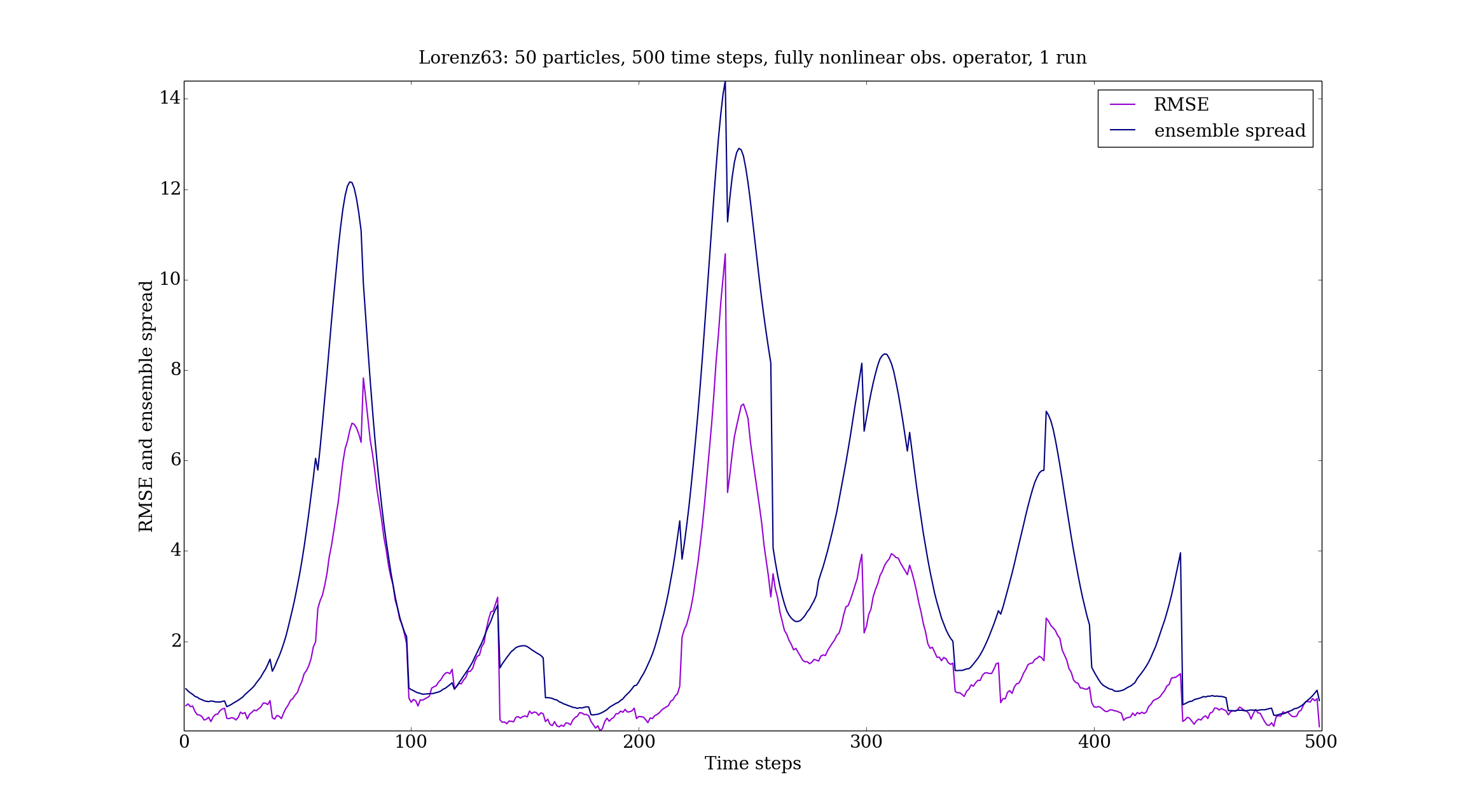}
   \caption{RMSE and ES}
   \label{daS9rmse}
  \end{subfigure}%
  \caption{Evolution of the Lorenz '63 model for 500 time steps, all 3 variables are observed every 20 time steps. Observations are denoted by $z_1,z_2,z_3$. Fully nonlinear observation operator i.e. $z_1=x^2, z_2=y^2, z_3=z^2$ plus noise.}
\end{figure}

\begin{figure}[ht!]
  \centering
  \begin{subfigure}{.5\linewidth}
    \centering
    \includegraphics[width = \linewidth]{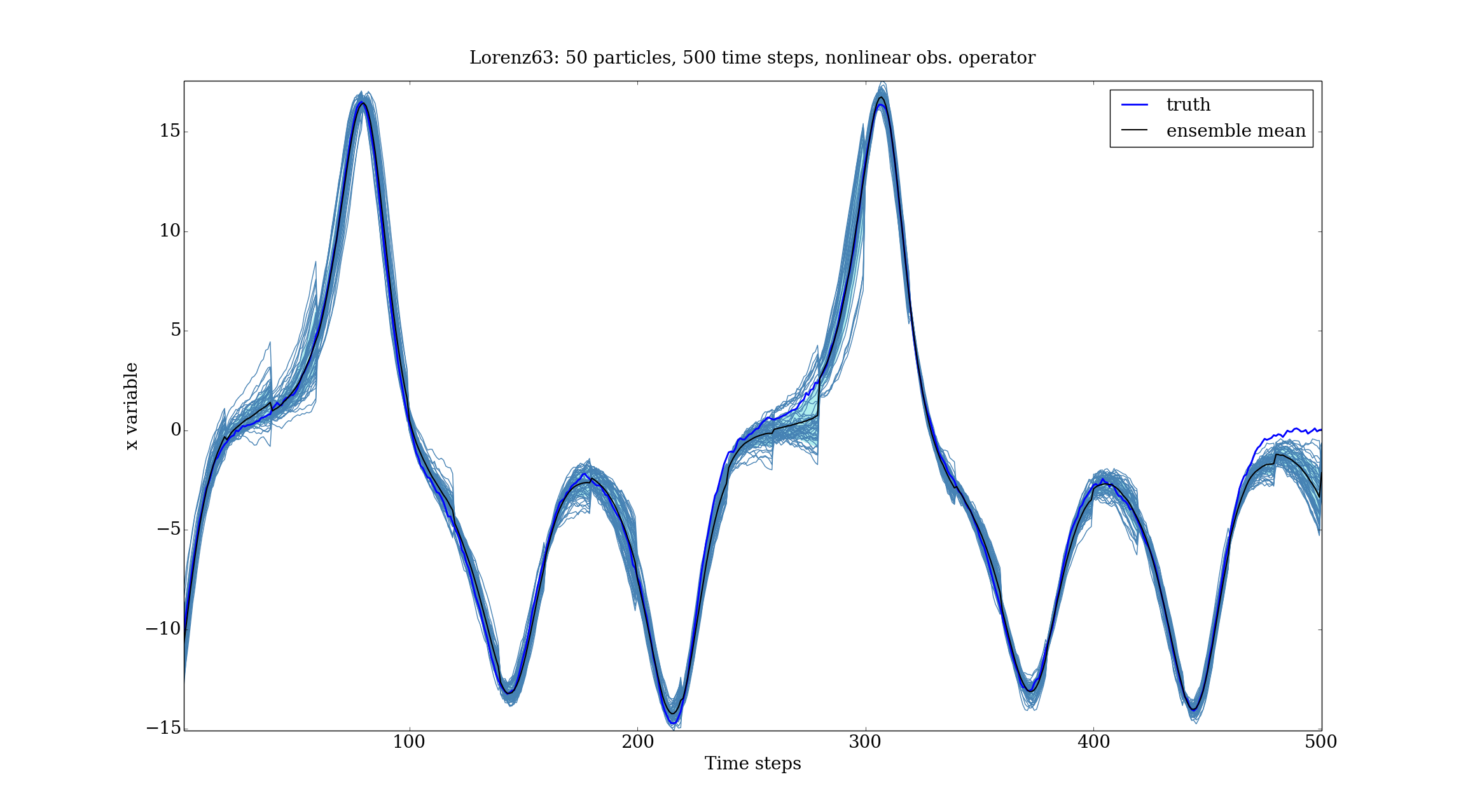}
    \caption{$x$ variable}
    \label{daS8updatedx}
  \end{subfigure}%
  \begin{subfigure}{.5\linewidth}
    \centering
   \includegraphics[width = \linewidth]{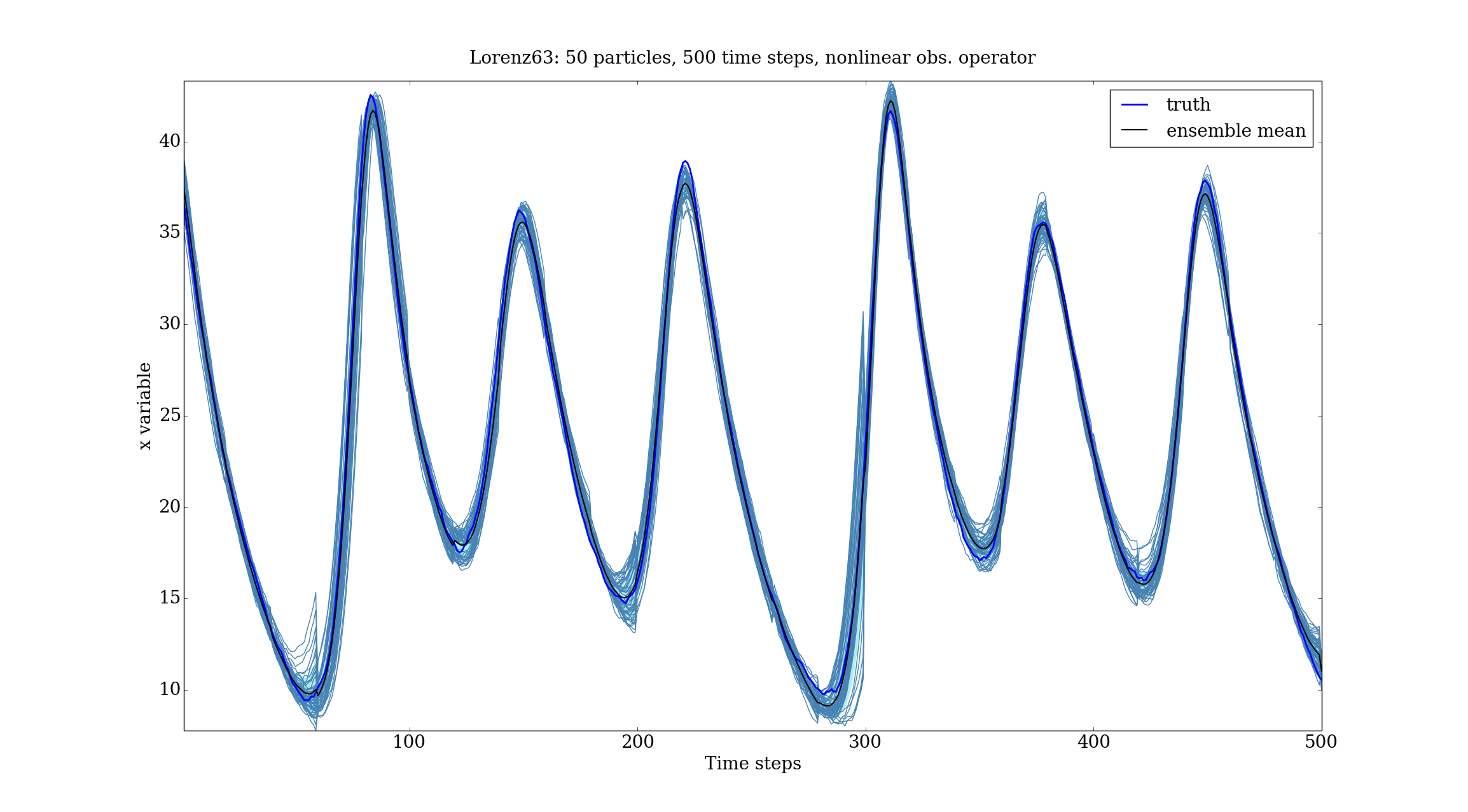}
   \caption{$z$ variable}
   \label{daS8updatedz}
  \end{subfigure}%
  \caption{Evolution of the Lorenz '63 model for 500 time steps, all 3 variables are observed every 20 time steps. Observations are denoted by $z_1,z_2,z_3$. Partially nonlinear observation operator i.e. $z_1=x^2, z_2=y, z_3=z$ plus noise. Displayed:  averaged (3 runs) RMSE and ensemble spread.}
\end{figure}

\begin{figure}[ht!]
  \centering
  \begin{subfigure}{.5\linewidth}
    \centering
    \includegraphics[width = \linewidth]{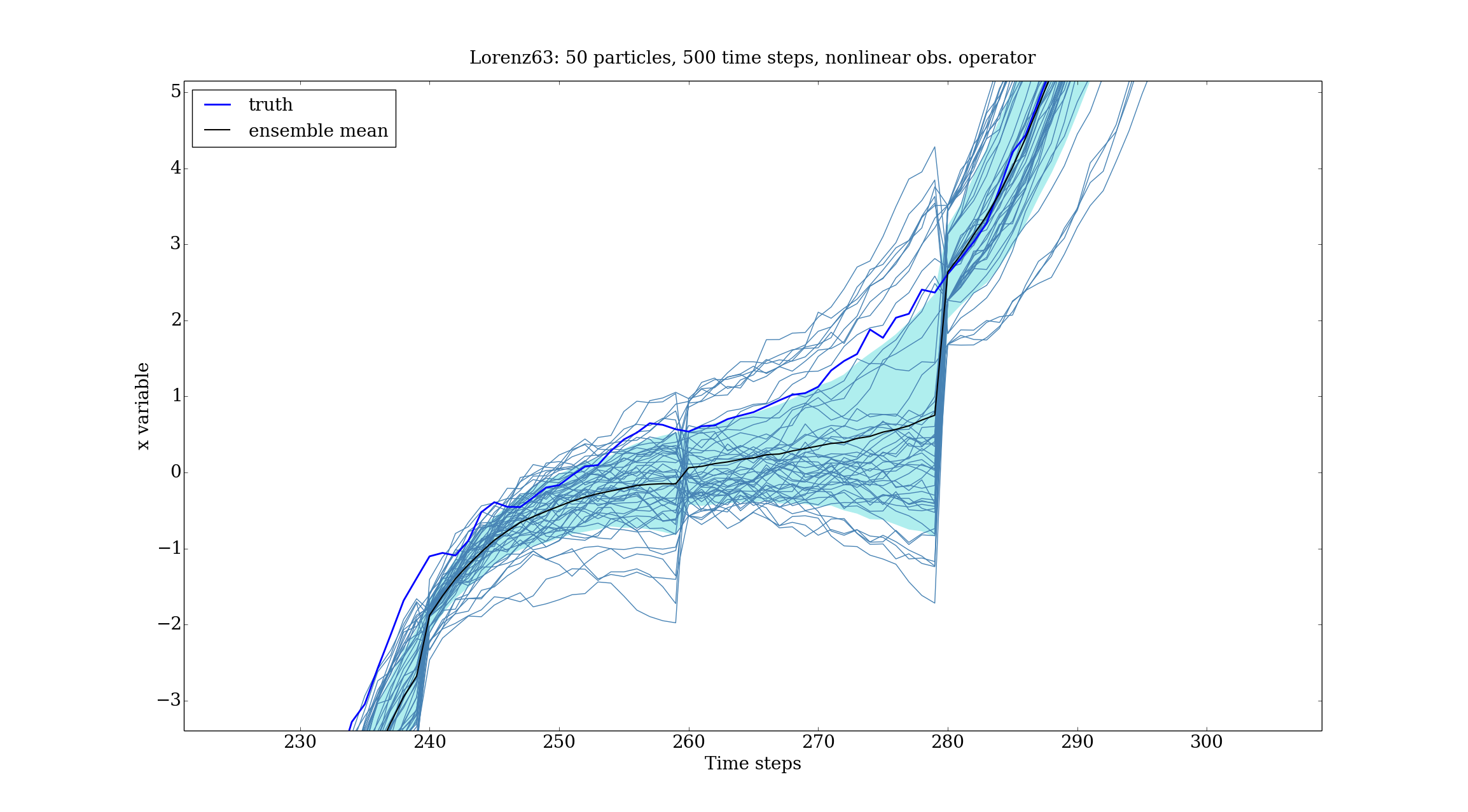}
    \caption{$x$ variable enlarged}
   \label{daS8updatedenlarged}
  \end{subfigure}%
  \begin{subfigure}{.5\linewidth}
    \centering
   \includegraphics[width = \linewidth]{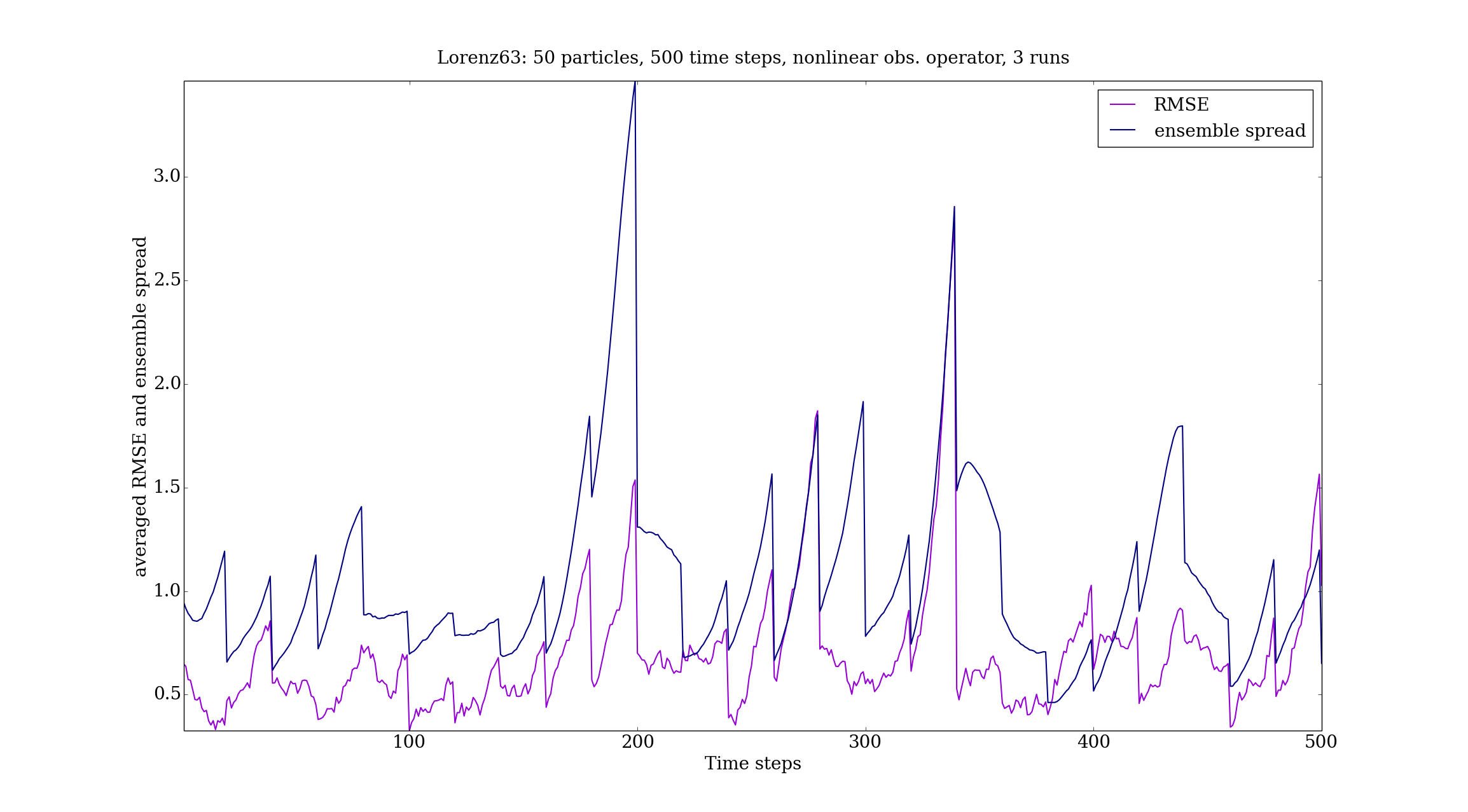}
   \caption{RMSE and ES}
   \label{daS8updatedrmse}
  \end{subfigure}%
  \caption{Evolution of the Lorenz '63 model for 500 time steps, all 3 variables are observed every 20 time steps. Observations are denoted by $z_1,z_2,z_3$. Partially nonlinear observation operator i.e. $z_1=x^2, z_2=y, z_3=z$ plus noise. Displayed:  averaged (3 runs) RMSE and ensemble spread.}
\end{figure}

In all previous settings the observation operator was linear. We now perform a couple of tests {(Figures \ref{daS9x}-\ref{daS8updatedrmse})} with nonlinear observation operators. We first make the observations fully nonlinear {(Figures \ref{daS9x}-\ref{daS9rmse})}. More precisely, the observation operator becomes
$$\mathscr{H}(x,y,z)=(x^2,y^2,z^2)$$
with observations given by
$$Z=\mathscr{H}(X_t)+V_t$$
as explained in Section \ref{particlefiltersection}. As we can see from Figure \ref{daS9x} the information available is not enough to keep the posterior concentrated around the truth at all times, as the system misses information in which wing of the butterfly it is. This is the case particularly with the $x$ direction. All particles follow the attractor, while some jump from one wing of the attractor to the other.
The behaviour of the posterior in the $z$ direction is much better, as the $z$ variable is always positive, and hence only one solution is present at all times for this variable. 
This can be observed in Figure \ref{daS9z}. The values of the RMSE and the ES are plotted in Figure \ref{daS9rmse}. We can see that both can become very large, essentially of the order of the support of the diffusion (the attractor), as the particles can spread around the entire attractor. 

Note that these results do not point to a failure of the particle filter. The exact posterior will also show this bimodal behavior, and the particles keep following the attractor as they should. Methods based on linearizations, such as Ensemble Kalman Filters, would show a cloud of particles between the two wings, which is the wrong solution.

To study this further we make 
the observation nonlinear just for the first variable {(Figures \ref{daS8updatedx}-\ref{daS8updatedrmse})}. That is, 
$$\mathscr{H}(x,y,z)=(x^2,y,z).$$
In contrast to the previous case, the particle filter performs much better. The (partial) linearity of the observation operator imposes a massive uncertainty reduction.
The cloud of particles remains concentrated around the signal trajectory and the data assimilation steps keep the uncertainty in check. Both the $x$ and the $z$ variable (Figures \ref{daS8updatedx}-\ref{daS8updatedz}) are successfully tracked, as highlighted in Figure \ref{daS8updatedenlarged}. The RMSE and the ES plotted in Figure \ref{daS8updatedrmse} confirm the findings: most of the time, the RMSE remains in the interval $[0, 1]$ with rare excursions away from this interval. The ensemble spread remains also small, with some oscillations around the assimilation time. The underlying reason for this behaviour is that in this case the system knows in which wing it is via the linear $y$ observation.


\section{Application to the Stochastic Rotating Shallow Water Model}
\subsection{Model Description}
The rotating shallow water model (RSW) is a classical nonlinear fluid dynamics model which contains key aspects of the oceanic and atmospheric dynamics. A detailed analytical description of this model has been provided in \cite{CL3}. From a numerical perspective, challenges are generated especially by the nonlinear advective terms. Nonlinear advection is a dispersive process. While in a linear case the wave-like solutions have constant amplitude and propagate at constant speed, in the nonlinear setting waves of different wavenumbers can propagate at varying speeds. The short waves can be amplified by the nonlinear structure, with direct impact on the accuracy of the solution (\cite{Durran}). These intricacies can be overcome by making use of the natural dominant balances which appear in atmosphere and oceans (geostrophic and hydrostatic balance) and by implementing the models using suitable numerical schemes. 



The stochastic rotating shallow water model used in our work is given by the following set of equations: 
\begin{subequations} \label{model}
 \begin{equation}\label{model1} 
 dv_t + \big[u_t \cdot \nabla v_t + f\hat{z} \times u_t + \nabla p_t \big]dt + \displaystyle\sum_{i=1}^{\infty}\big[(\mathcal{L}_i + \mathcal{A}_i)v_t\big] \circ dW_t^i = 0 
 \end{equation}
 \begin{equation}\label{model2} 
  dh_t + \nabla \cdot (h_t u_t)dt+\displaystyle\sum_{i=1}^{\infty}\big[\nabla \cdot (\xi_i h_t)\big]\circ dW_t^i= 0 
 \end{equation}
\end{subequations}
where 
$\mathcal{L}_iv := \xi_i \cdot \nabla v, \ \mathcal{A}_iv:= v_j\nabla \xi_i^j = \displaystyle\sum_{j=1}^2 v_j\nabla\xi_i^j$, $v=(v^1,v^2)$,$v: = \epsilon u + \mathcal{R}$,
$u=(u^1,u^2)$ is the horizontal fluid velocity vector,
   $p:=\frac{h-b}{\epsilon \mathscr{F}}$ is the pressure term, $h$ is the total depth, $b=b(x)$ is the bottom topography, $x=(x^1,x^2)$, $f=f(x)$ is the Coriolis parameter, $\hat{z}$ is a unit vector pointing away from the centre of the Earth, $\epsilon <<1$ is the Rossby number, $\mathscr{F}$ is the Froude number, $\mathcal{R}=\mathcal{R}(x)$ is the vector potential of the divergence-free rotation rate about the vertical direction, with $curl \ \mathcal{R}(x)=f(x)\hat{z}.$
   
{ We start with a velocity vector field which is in geostrophic balance at the initial time $t=0$. Nonetheless, we eliminate the geostrophic balance condition on the velocity field from $t=1$ onwards. On the other hand, the random fields are assumed to stay in geostrophic balance at each time step and we give the explicit formula for this in the next subsection, equations \eqref{balance}}.
Intuitively, the geostrophic balance can be explained as follows (\cite{VallisOceans} pp. 58): initially there exists a pressure gradient within the fluid, which generates a fluid dynamics from high-pressure regions towards low-pressure regions; while the fluid flow evolves in time, it is deflected by the Coriolis force; in the Northern hemisphere the Coriolis force deflects the fluid to the right, while the pressure gradient does not change direction much, so that the two forces can come into equilibrium. The resulting direction of the fluid motion is perpendicular on both the Coriolis force and the pressure force.
  In the atmosphere and oceans this balance tends to be stable, especially away from the boundaries. Small disturbances to it lead to the appearance of gravity waves. The geostrophic balance is dominant as the wind and the currents are usually weak in comparison to the speed of the Earth rotation (\cite{VallisOceans}, pp.59). The Rossby number in this case is small ($\epsilon \sim 10^{-1}$), meaning that the rotation dominates the advective part and it is balanced by the pressure gradient force (\cite{MPEbook}, pp. 86). If the rotation is fast, then the horizontal flow is in 'near' geostrophic  balance, that is 'nearly' divergence-free (\cite{Vallis}, pp.95). We also assume geostrophic balance for the stochastic forcing described in Section \ref{stochasticforcing}, to strongly reduce the generation of artificial gravity waves.  

{ We initialise the model by first calculating the Coriolis parameter and the pressure term for each grid point. Then we generate $v_0=(v_0^1,v_0^2)$ using the geostrophic balance condition.} The periodicity of the domain inspires a rigorous numerical setup for generating the stochastic forcing. This is described in Section \ref{stochasticforcing}. The system \eqref{model} is implemented on a staggered grid, using a Runge-Kutta scheme of order 4. The domain corresponds to a strip situated between 30 and 60 degrees north latitude. The numerical scheme is described in Appendix. A staggered grid (also known as Arakawa C-grid) is usually preferred when implementing weather prediction models due to the low dispersion errors and the lack of computational modes associated with it (\cite{MPEbook}, pp. 313).  

\subsubsection{Stochastic forcing} \label{stochasticforcing}
   For a state vector $X$ of dimension $d$ the nonlinear stochastic forcing is given by
   \begin{equation}
      \displaystyle\sum_{i=1}^{\infty}\mathcal{B}_i(\xi_i)X_tdW_t^i \sim \sqrt{\mathscr{Q}}dW (X) \sim \left(\displaystyle\sum_{i=1}^{\infty}\mathcal{B}_i(\xi_i)dW^i\right)X =: R_i X
      \end{equation}
  where $\mathcal{B}_i:L^2(\T^2) \rightarrow L^2(\T^2)$ is a general operator which depends on the vector fields $\xi_i$. In our case $\mathcal{B}_i$ can be given by $\cL_i$ or by $\cL_i + \cA_i$. With $\sim$ we denoted the equivalence in notation of the terms above, it has no mathematical meaning other than expressing different ways of rewriting the same stochastic term. Then we interpret $\displaystyle\sum_{i=1}^{\infty}\mathcal{B}_i(\xi_i)dW_t^i$ as a random operator applied to the state vector $X$ and we denote it by $R_i$. Here
  $dW$ is a $d$-dimensional vector of independent Brownian motions and
   $\mathscr{Q}$ is a space-covariance operator. We discretize the SPDE in time and space, and therefore $\mathscr{Q}$ becomes a covariance matrix of dimension $d^2$.

  \noindent Given the fact that generating $\sqrt{\mathscr{Q}}dW$ is computationally expensive, we first do this in the spectral space, and then return to the physical space. The covariance matrix is symmetric and circulant. We determine it in the spectral space for a system which is periodic in both the $x$ and $y$ directions, using the fast Fourier transform. Since our original domain is not periodic, we then truncate the resulting field to the physical domain, in order to avoid a periodic random field. In particular, we compute the Fourier transform of a column corresponding to the circulant matrix. This column has a Gaussian correlation structure. At every time step, we generate a Gaussian random field (again in the spectral space) and then perform the multiplication between the column of the covariance matrix and this newly generated random field. Finally, we use the inverse fast Fourier transform to return to the physical domain. 
  
  It is known (see \cite{RuanMcLaughlin}) that any continuous random field $R_i$ can be expressed as a Fourier-Stieltjes integral over a complex-valued Fourier increment $dY_{R_i}$:
   \begin{equation}\label{Ri}
   R_i(\bx) = R_i(x,y) = \displaystyle\int_{-\infty}^{\infty} e^{i\bmu\cdot \bx}dY_{R_i}(\bmu) \quad i=1,2, \ldots
   \end{equation}
   where the random component $Y_{R_i}$ must satisfy the following properties:
   \begin{equation}
   \begin{aligned}
   &\mathbb{E}[Y_{R_i}(\bmu)]=0 \\
   & \mathbb{E}[Y_{R_i}(\bmu)Y_{R_j}^{*}(\bmu')]=0 \quad \hbox{for} \quad \bmu \neq \bmu'\\
   & \mathbb{E}[Y_{R_i}(\bmu)Y_{R_j}^{*}(\bmu')]=Q_{ij}(\bmu)d\bmu
   \end{aligned}
   \end{equation}
   where $Y_{R_j}^{*}$ is the complex conjugate of $Y_{R_j}$, $Q_{ij}$ is the covariance between the two processes $R_i$ and $R_j$, $n_{\bx}$ is the dimension of the spatial vector $\bx$, and $\bmu = (\mu_1, \mu_2, \ldots, \mu_{n_{\bx}})$ corresponds to the wave number. Therefore, the key ingredient in generating the random fields $(R_i)_i$ is the accurate retrieval of the processes $Y_{R_i}(\bmu)$. In practice this should be effectuated in discrete time, so that one can implement the integral $\eqref{Ri}$ numerically. Inspired by the periodic structure of the domain, one way of efficiently performing this is by using a complex-valued spectral decomposition (\cite{Evensen}). We know from classical Fourier analysis that the vector fields $R_i$ can also be expressed as
   \begin{equation}
   R_i(\bx) = \displaystyle\int_{-\infty}^{\infty} e^{i\bmu \cdot \bx}\hat{R}_i(\bmu)d\bmu
   \end{equation}
   where $\hat{R}_i$ is the Fourier transform of $R_i$.
    Therefore we can write
      \begin{equation}\label{stochrandomfield}
      R_i(x_{k_1}, y_{k_2}) = \displaystyle\sum_{j,k}\hat{R}_i(\mu_j,\mu_k) e^{i(\mu_jx_{k_1}+\mu_ky_{k_2})} \Delta \bmu
      \end{equation} 
      with $x_{k_1} = k_1\Delta x$, $y_{k_2}=k_2\Delta y$, $\mu_j = \frac{2\pi j}{n\Delta x}$, $\mu_k = \frac{2\pi k}{m\Delta y}$, $\Delta \mu = \Delta \mu_j \Delta \mu_k  =\frac{(2\pi)^2}{nm \Delta x\Delta y}$ and
      \begin{equation}
       \hat{{R}_i}(\mu_j , \mu_k) = \frac{C}{\sqrt{\Delta\bmu}}e^{-\frac{\mu_j^2+\mu_k^2}{\sigma^2}+2\pi i \varphi_{jk}}
      \end{equation} 
      where $(n,m)$ is the dimension of the grid and $\varphi_{jk} \in [0,1]$ is a random number which introduces a random phase shift for any given wave number. 
      Therefore
      \begin{equation}
      \begin{aligned} 
      R_i(x_{k_1}, y_{k_2}) &= C \sqrt{\Delta \bmu} \displaystyle\sum_{j,k}e^{-\frac{\mu_j^2+\mu_k^2}{\sigma^2}+2\pi i \varphi_{jk}} e^{i(\mu_jx_{k_1}+\mu_ky_{k_2})} \\
      & = |\Delta\bmu|^{1/2}\displaystyle\sum_{j,k}e^{\varphi_{jk}(\bmu)}C e^{i\bmu \cdot \bx}\\
      & \approx \Delta Y_{R_i}(\bmu) \approx dY_{R_i}(\bmu)
      \end{aligned} 
      \end{equation}
  where $\varphi_{jk}$ is a random process which corresponds to the (discretized) wave number domain with volume element given by $|\Delta\bmu| = \displaystyle\prod_{j=1}^{n_x}\Delta \bmu_j$ for any fixed wave number $\bmu$. Therefore
  \begin{equation}\label{Rifinal}
  \begin{aligned}
  R_i(\bx) &= \displaystyle\int_{-\infty}^{\infty} e^{i\bmu \cdot \bx} dY_{R_i}(\bmu)\\
  & \approx \displaystyle\sum_{j,k}e^{\varphi_{jk}(\bmu)}Ce^{i\bmu\cdot\bx} |\Delta \bmu|^{1/2}
  \end{aligned}
  \end{equation}
  that is, a summation over all the discrete wave numbers of the grid. By the central limit theorem, $R_i$ are normally distributed. In order for the expression in \ref{Rifinal} to be well-defined, $\Delta Y_{R_i}$ and $\Delta Y_{R_j}$ must be orthogonal in the limit as the discretization becomes finer (\cite{RuanMcLaughlin}). An efficient way of computing \ref{Rifinal} is by using a fast Fourier transform algorithm. One has
  \begin{equation}
  \begin{aligned} 
  \mathbb{E}[R_i(x)R_j^{*}(x')] &= \displaystyle\int_{-\infty}^{\infty}\displaystyle\int_{-\infty}^{\infty} e^{i\bmu \cdot \bx - i\bmu' \cdot \bx'}\mathbb{E}[dY_{R_i}(\bmu)dY_{R_j}^{*}(\bmu')]\\
  & = \displaystyle\int_{-\infty}^{\infty}e^{i\bmu(\bx-\bx')}Q_{ij}(\bmu)(d\bmu)
  \end{aligned} 
  \end{equation} 
  that is the covariance of $R_i$ and $R_j^{*}$ can be written as the inverse Fourier transform of $Q_{ij}$. 
  More details on random field generators can be found in \cite{RuanMcLaughlin} and \cite{Evensen}.  
  { As mentioned in the previous section,} we assume that the dominant balance { for the random vectors} in this system is the \textit{geostrophic balance}, that is the pressure gradient force and the Coriolis force balance each other. { We have three random vectors: $R^{v^1}, R^{v^2}$ corresponding to each component of the velocity $v=(v^1,v^2)$, and $R^p$ which corresponds to the pressure field $p$. We obtain them as follows: first, we compute $R^p$ using equation \eqref{stochrandomfield}; then we use $R^p$ to compute $R^{v^1}$ and $R^{v^2}$ under the geostrophic balance assumption:}
  {\begin{subequations}\label{balance} 
  \begin{equation}\label{balance1} 
   R^{v^1} = -\frac{g}{f}\frac{\partial R^p}{\partial y} 
  \end{equation}
  \begin{equation}\label{balance2}
   R^{v^2} = \frac{g}{f}\frac{\partial R^p}{\partial x}.
  \end{equation}
  \end{subequations}
  } 
\subsubsection{Connection with SALT}
The SRSW model considered here has been derived using the SALT approach described in \cite{Holm2015} (see Appendix). Stochasticity is introduced in the advection part of the dynamics to model the uncertain transport behaviour within the fluid flow. To the best of our knowledge this is the first implementation of the SALT SRSW model in a data assimilation setting. Numerical implementations and particle filter algorithms for other SALT models (2D Euler, SQG) have been extensively developed in \cite{Wei1}-\cite{Weifinal}.  
Our numerical implementation of the stochastic forcing is equivalent to the one described in \cite{Wei1} up to an isomorphism. More precisely, if one extends the approach described in \cite{Wei1} to a periodic domain, then the stochastic vector fields correspond to the Fourier modes of the transport covariance matrix and they form a basis of the underlying space. This basis may be different from the basis chosen by us above. However, an isomorphism can be established between any two orthogonal bases corresponding to the same underlying space. \\

\subsection{Data Assimilation Results}
\begin{figure}[ht!]
  \centering
  \begin{subfigure}{.5\linewidth}
    \centering
    \includegraphics[width = \linewidth]{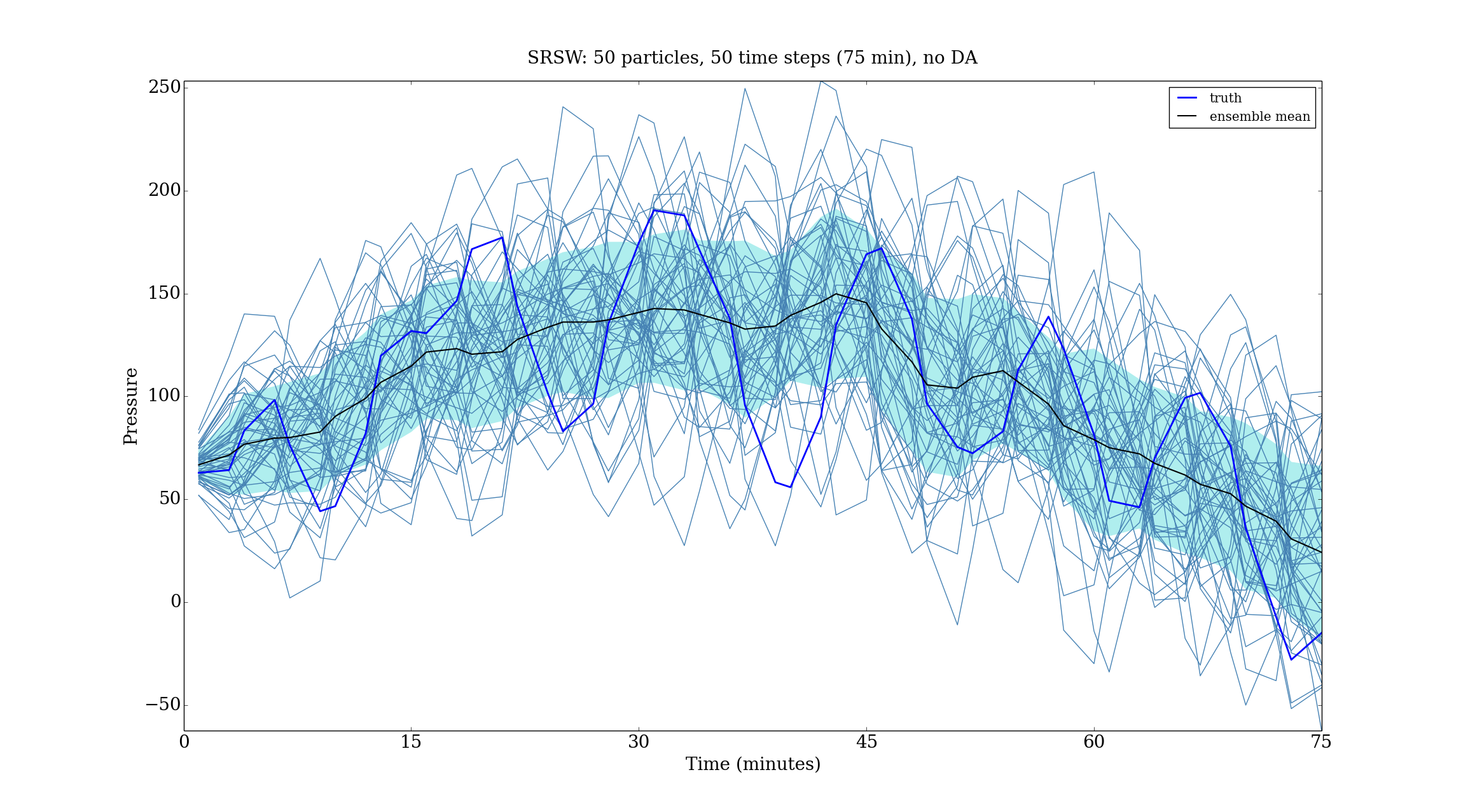}
    \caption{pressure field}
    \label{noDAsrswS1-1p}
  \end{subfigure}%
  \begin{subfigure}{.5\linewidth}
    \centering
   \includegraphics[width = \linewidth]{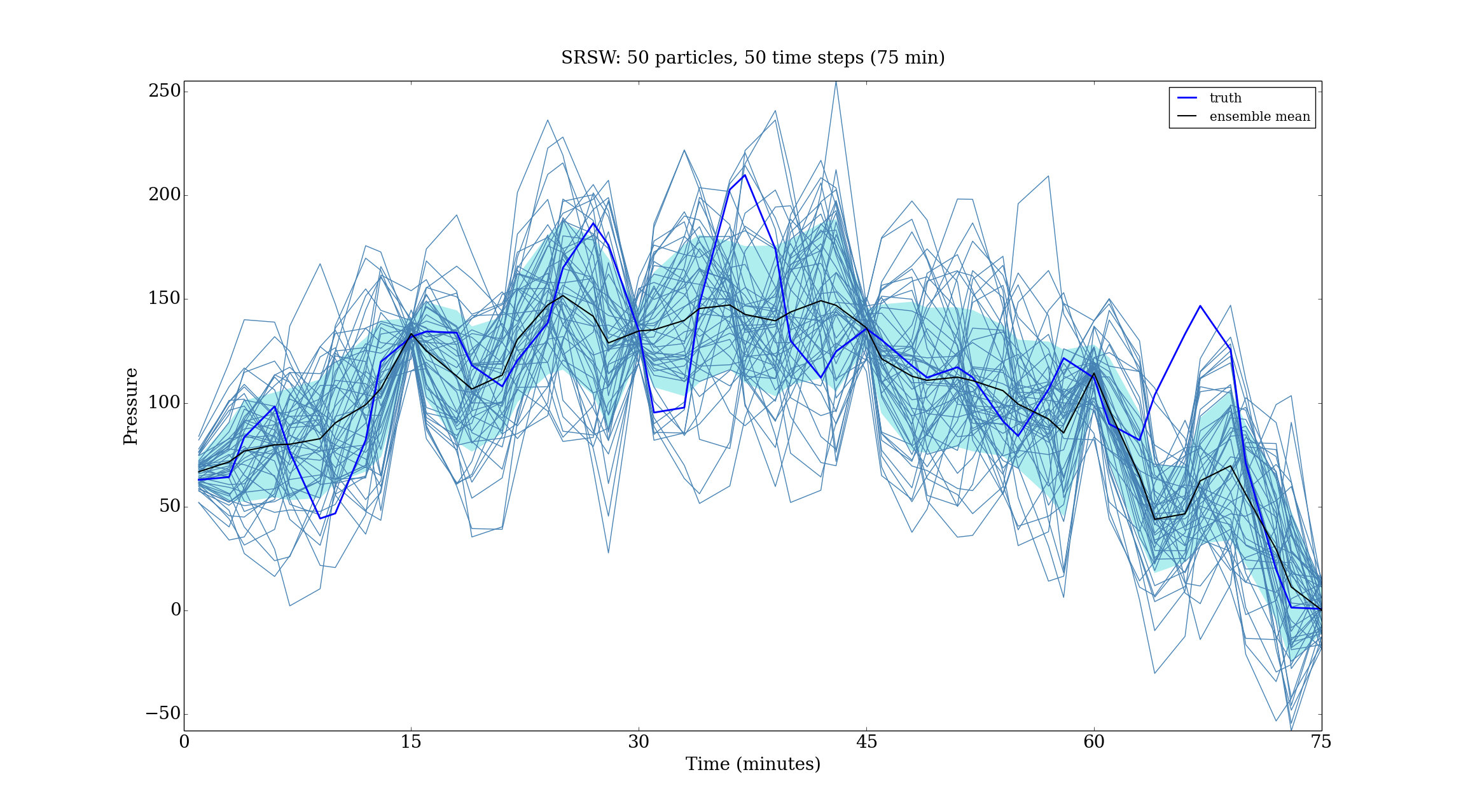}
   \caption{pressure field}
  \label{dasrswS1p}
  \end{subfigure}%
\end{figure}

\begin{figure}[ht!]
  \centering
  \begin{subfigure}{.5\linewidth}
    \centering
    \includegraphics[width = \linewidth]{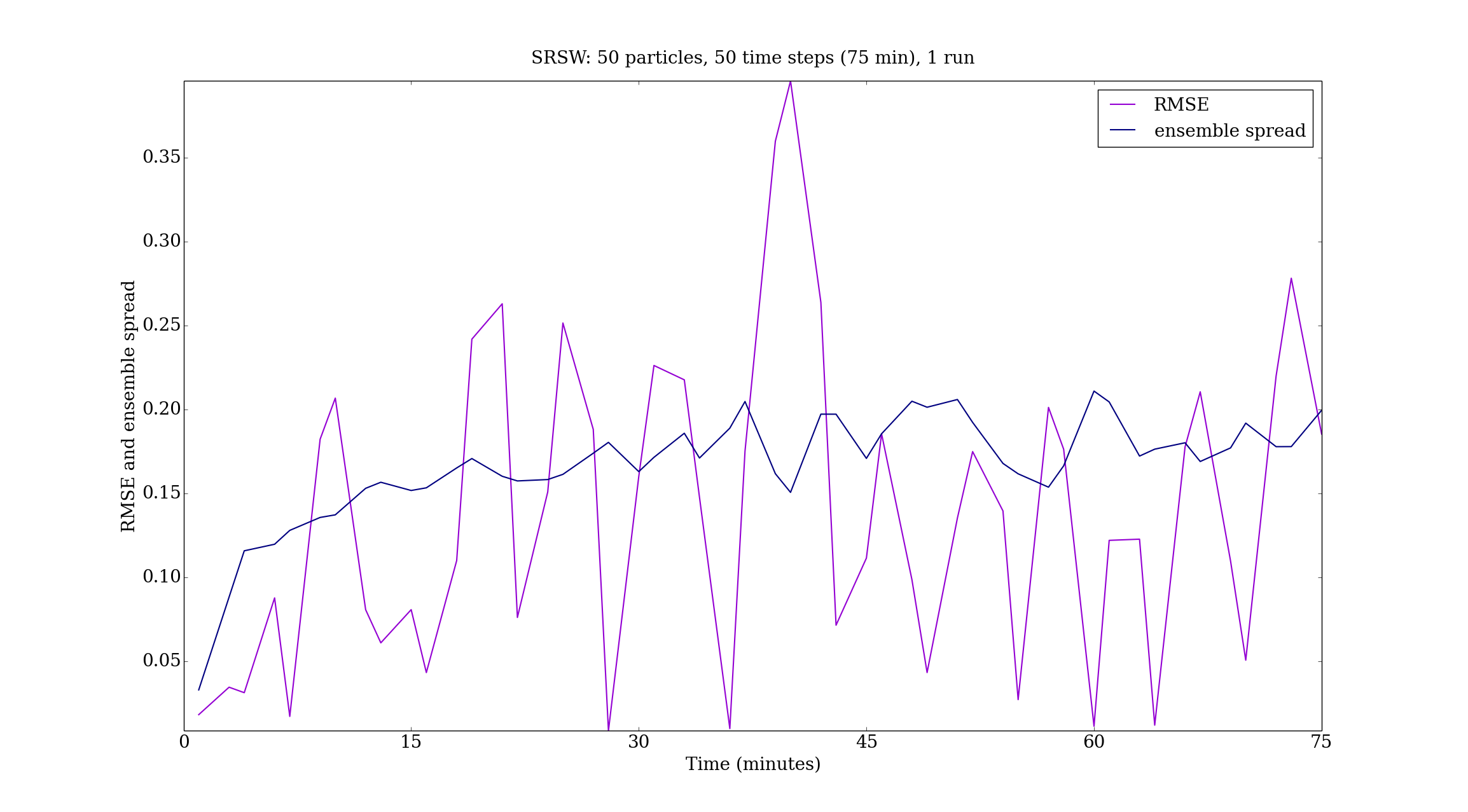}
    \caption{RMSE}
   \label{noDAsrswS1-1rmse}
  \end{subfigure}%
 \begin{subfigure}{.5\linewidth}
    \centering
   \includegraphics[width = \linewidth]{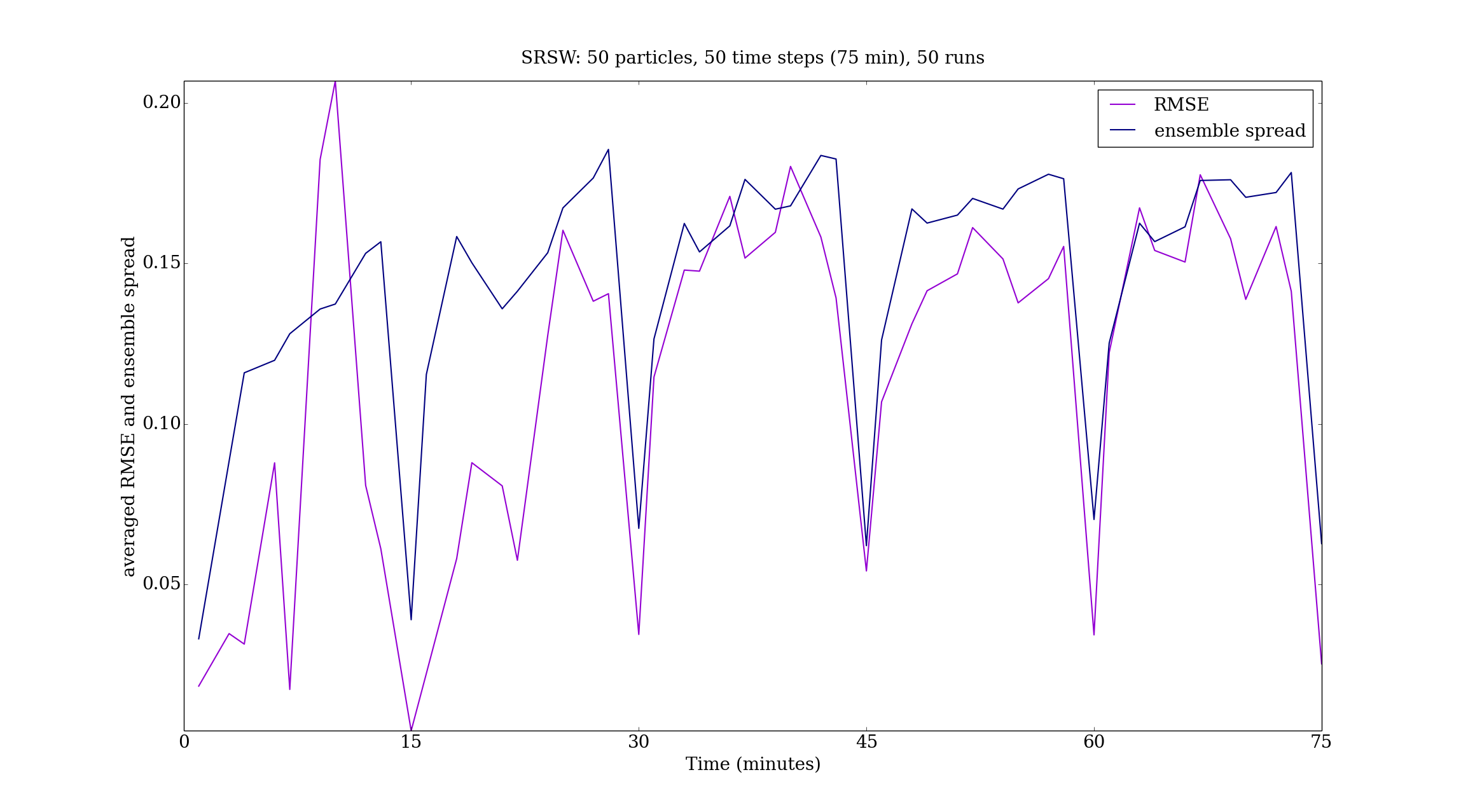}
   \caption{RMSE}
   \label{dasrswS1rmse}
  \end{subfigure}%
  \caption{Evolution of the SRSW model for 50 time steps. Left: no assimilation of data. Right: the system is observed every 10 time steps.}
\end{figure}

We perform the data assimilation analysis using an ensemble of 50 particles. 
We plot the ensemble of particles one standard deviation region about the ensemble mean and compare it with the truth. The \textit{truth} is a pathwise realisation of the SRSW model. In the \textit{standard} setting {(Figures \ref{dasrswS1p}-\ref{dasrswS1rmse})} we observe the system for 50 time steps, with a time step size of 90 seconds (hence we observe every 1.25 hour). We use one observation at each analysis time, and the observational uncertainty and the initial uncertainty are equal to 1. The stochastic model error is set to 200 metres. As mentioned, the random forcing is implemented using a stochastic advection velocity in SALT. We generate these stochastic advection fields by first generating a stochastic pressure field of amplitude 200 m, then determining the velocity fields in the zonal and meridional directions, and using the geostrophic balance, as explained in the previous section.
 These then form the stochastic advection fields that are applied in all three evolution equations. 

We first show the output obtained when no data is assimilated. Figure \ref{noDAsrswS1-1p} contains the evolution of the ensemble of particles corresponding to the pressure field at a grid point in the middle of the domain. We can see that the prior distribution is spread out and the truth is not successfully tracked (the ensemble mean remains far from the truth).


\begin{figure}[ht!]
  \centering
  \begin{subfigure}{.5\linewidth}
    \centering
    \includegraphics[width = \linewidth]{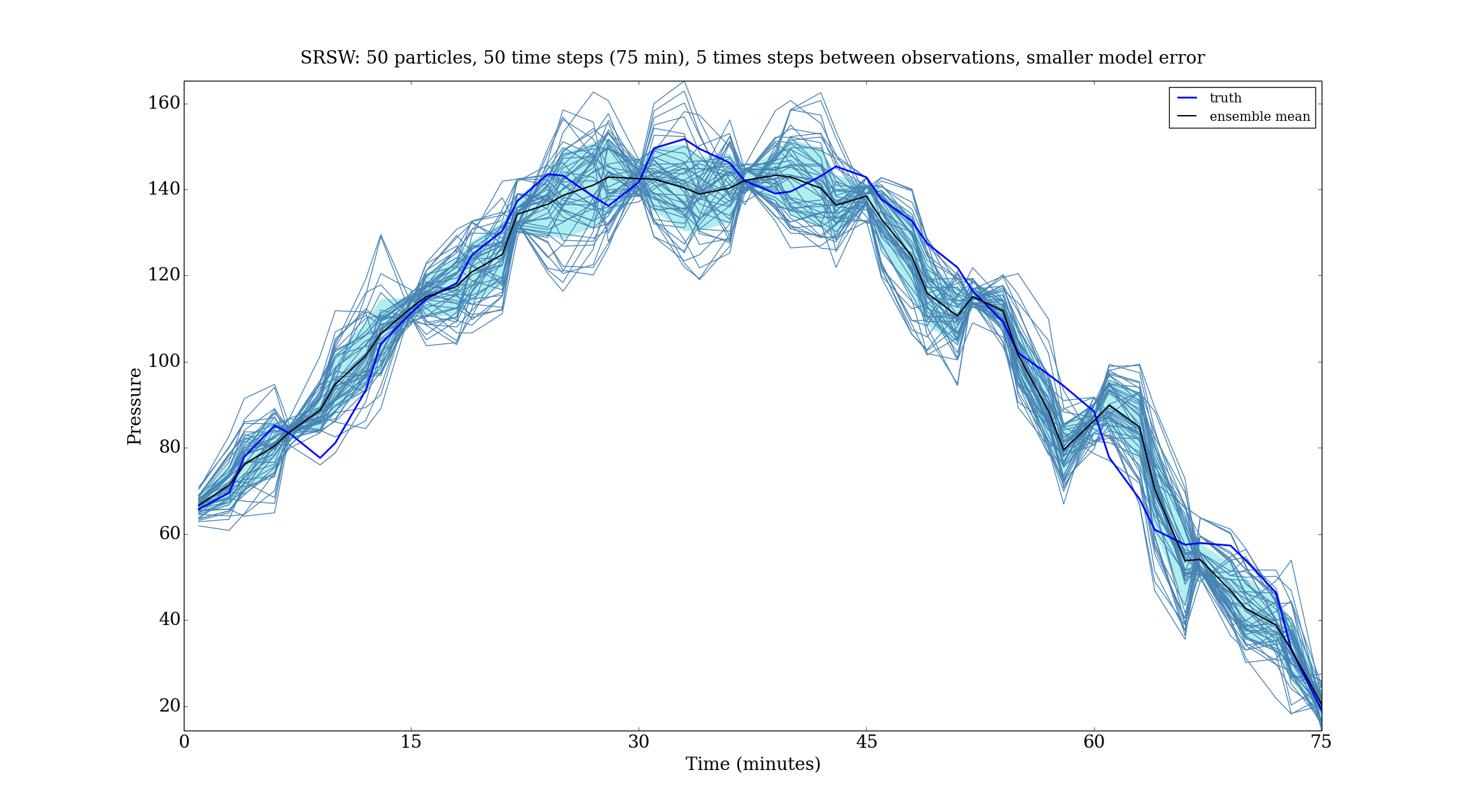}
    \caption{pressure field}
   \label{dasrswS4p}
  \end{subfigure}%
 \begin{subfigure}{.5\linewidth}
    \centering
   \includegraphics[width = \linewidth]{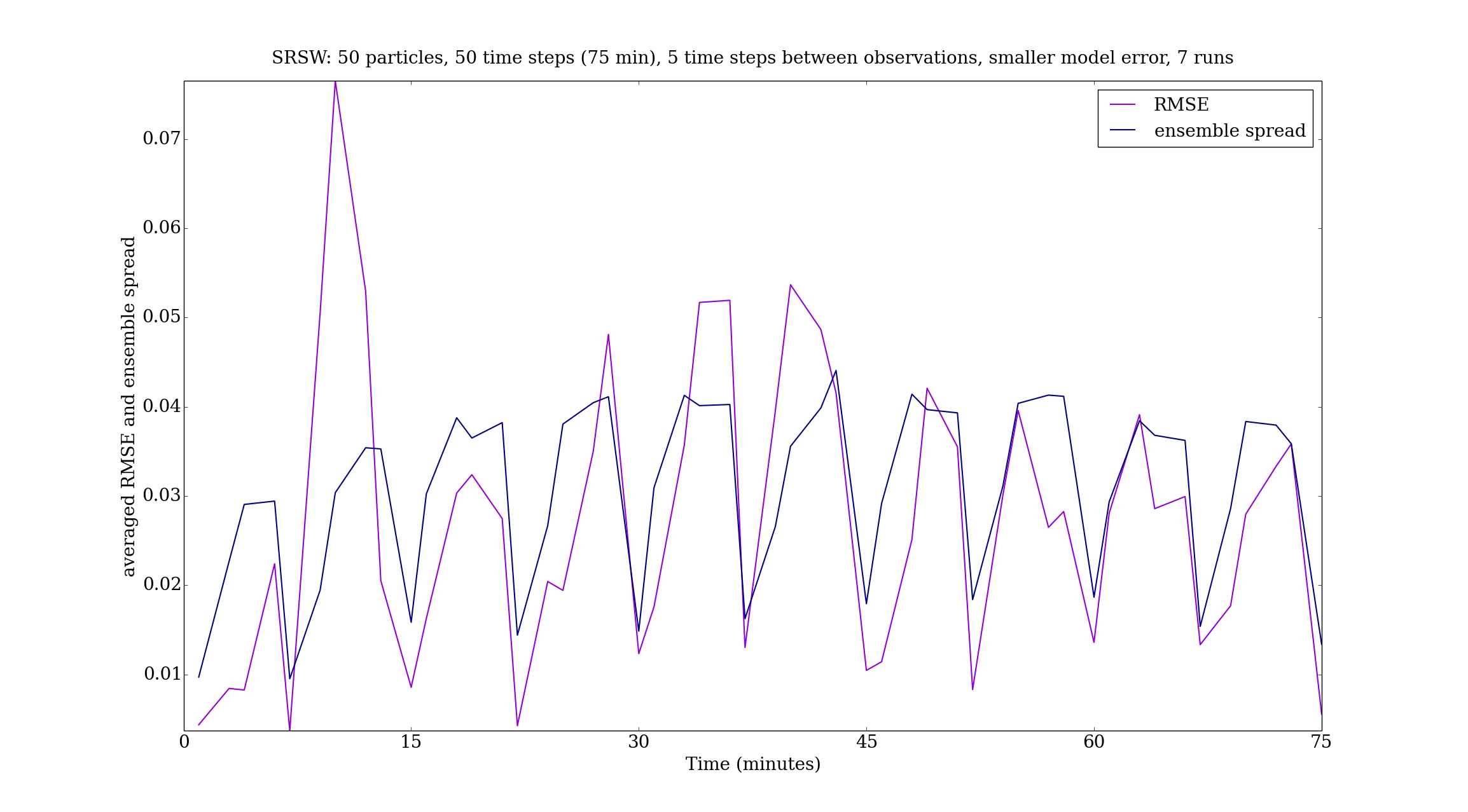}
   \caption{RMSE}
   \label{dasrswS4rmse}
  \end{subfigure}%
  \caption{Evolution of the SRSW model for 50 time steps. The system is observed every 5 time steps. The model error is decreased to 50. { The particle filter efficacy is improved by the time window and the decreased model error}.}
\end{figure}

Now we start to observe the system every 10 time steps {(Figures \ref{dasrswS1p}-\ref{dasrswS1rmse})}, so every 15 minutes. In Figure \ref{dasrswS1p} one can see that the assimilation of data improves the performance of the ensemble of particles: the standard deviation is reduced and the particle trajectories are \textit{corrected} at assimilation times. 
We then decrease the stochastic model error to 50 and observe the system every 5 time steps { (Figures \ref{dasrswS4p}-\ref{dasrswS4rmse})}. In figure \ref{dasrswS4p} one can notice that the particle filter efficiency is significantly improved by the data assimilation time window: by observing the system more often we obtain a cloud of particles which are better concentrated around the truth. This success is due also to the fact that we decreased the model error to 50. The RMSE and SE illustrated in Figure \ref{dasrswS4rmse} are much smaller than before.

\begin{figure}[ht!]
  \centering
  \begin{subfigure}{.5\linewidth}
    \centering
    \includegraphics[width = \linewidth]{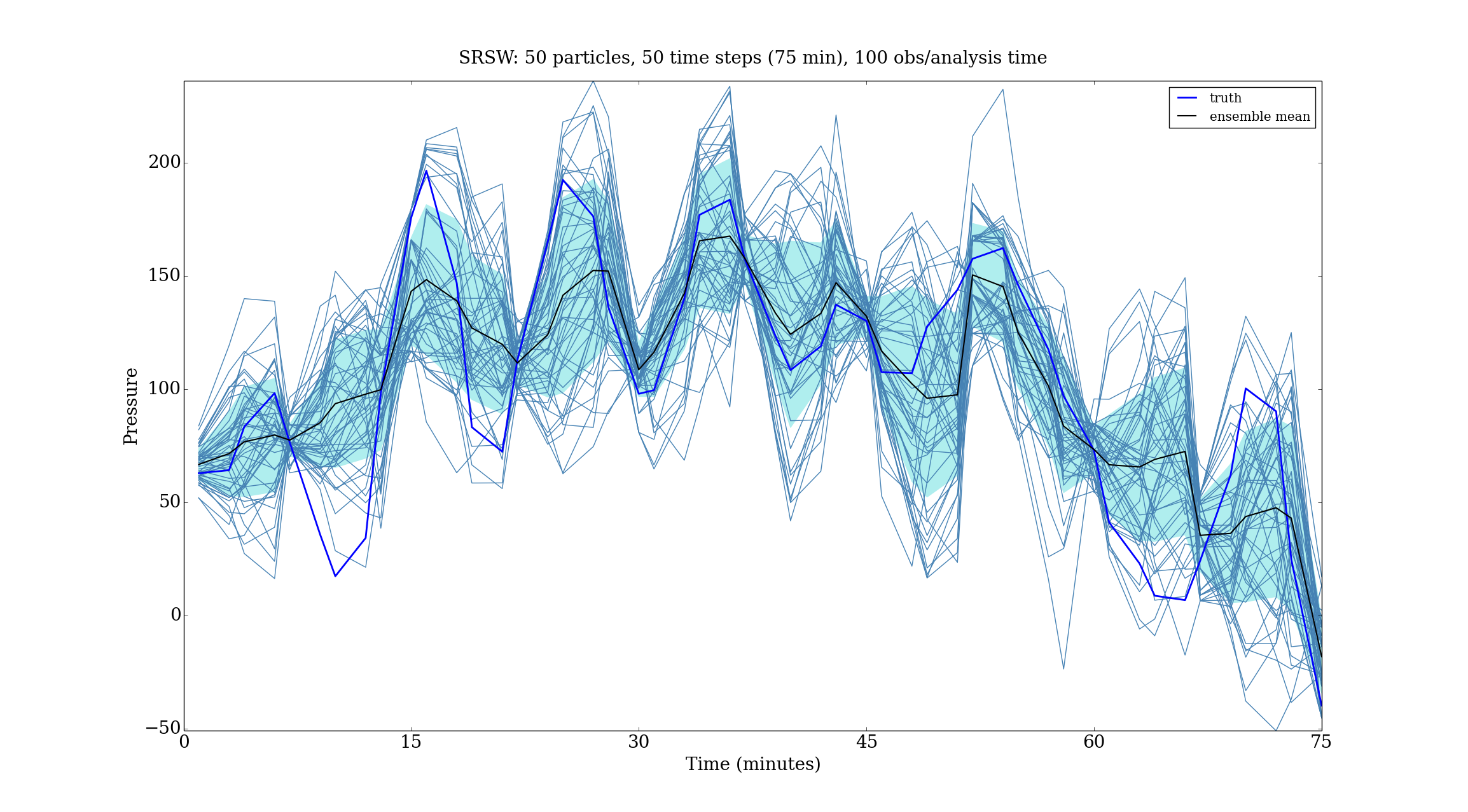}
    \caption{pressure field}
   \label{dasrswS6p}
  \end{subfigure}%
 \begin{subfigure}{.5\linewidth}
    \centering
   \includegraphics[width = \linewidth]{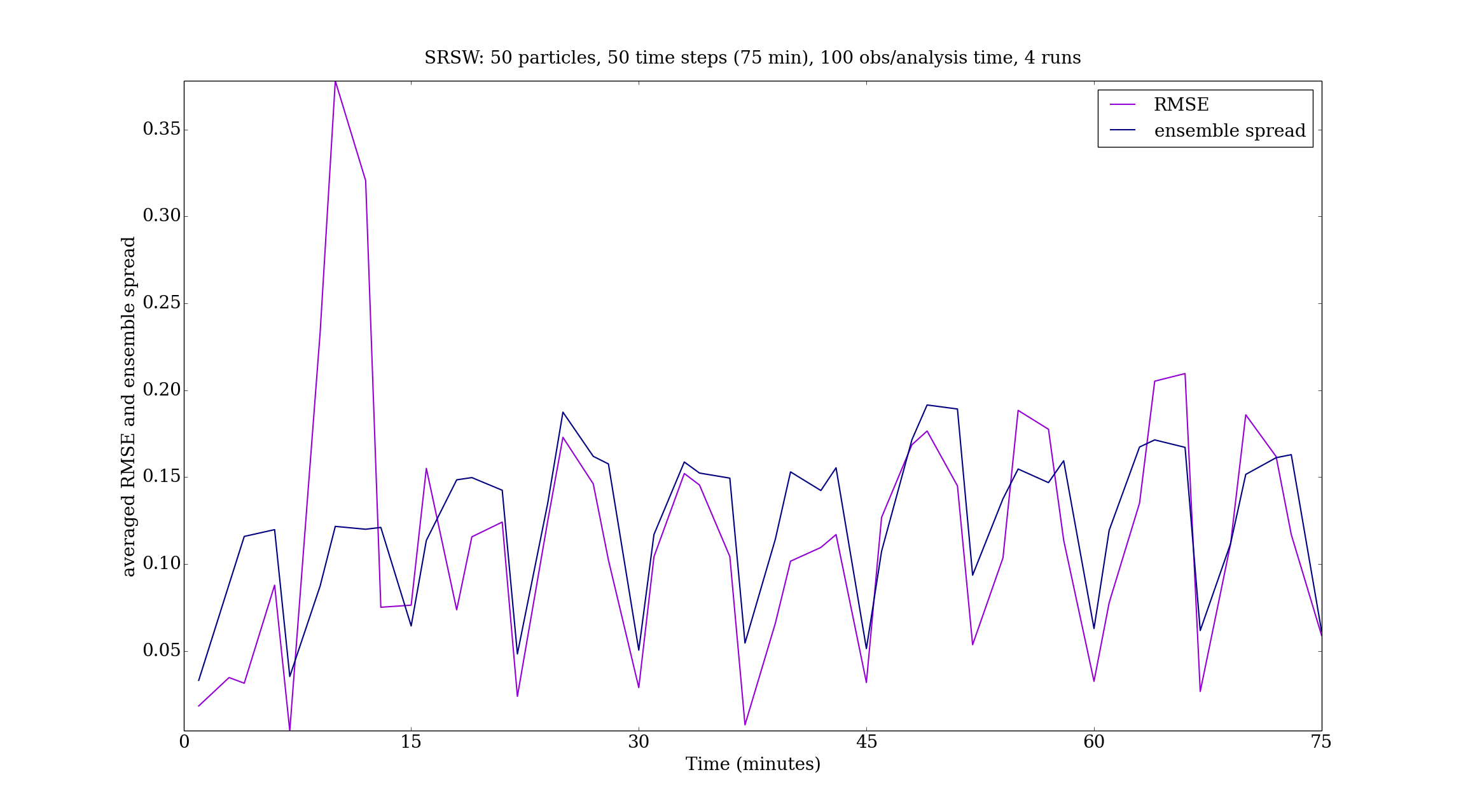}
   \caption{RMSE}
   \label{dasrswS6rmse}
  \end{subfigure}%
  \caption{Evolution of the SRSW model for 50 time steps. The system is observed every 5 time steps, with 100 observations/analysis time. { A large number of observations improves the performance of the particle filter.}}
\end{figure}

\begin{remark}
Note that all plots are rescaled automatically. 
\end{remark}
For the next scenario {(Figures \ref{dasrswS6p}-\ref{dasrswS6rmse})} the system is observed every 5 time steps and we use 100 observations at each analysis time. Figure \ref{dasrswS6p} shows that a drastic increase in the number of observations (100 now as opposed to 1 as we had before) greatly improves the efficiency of the particle filter. In Figure \ref{dasrswS6rmse} one can see that the ensemble spread is substantially reduced and the accuracy is now also much better.
The true innovation is that although we do not perform any localisation, we can assimilate 100 observations without filter degeneracy. This suggest that this particle filter is an important step towards beating the \textit{curse of dimensionality} that plagued particle filtering for so long. 

\begin{figure}[ht!]
  \centering
  \begin{subfigure}{.5\linewidth}
    \centering
    \includegraphics[width = \linewidth]{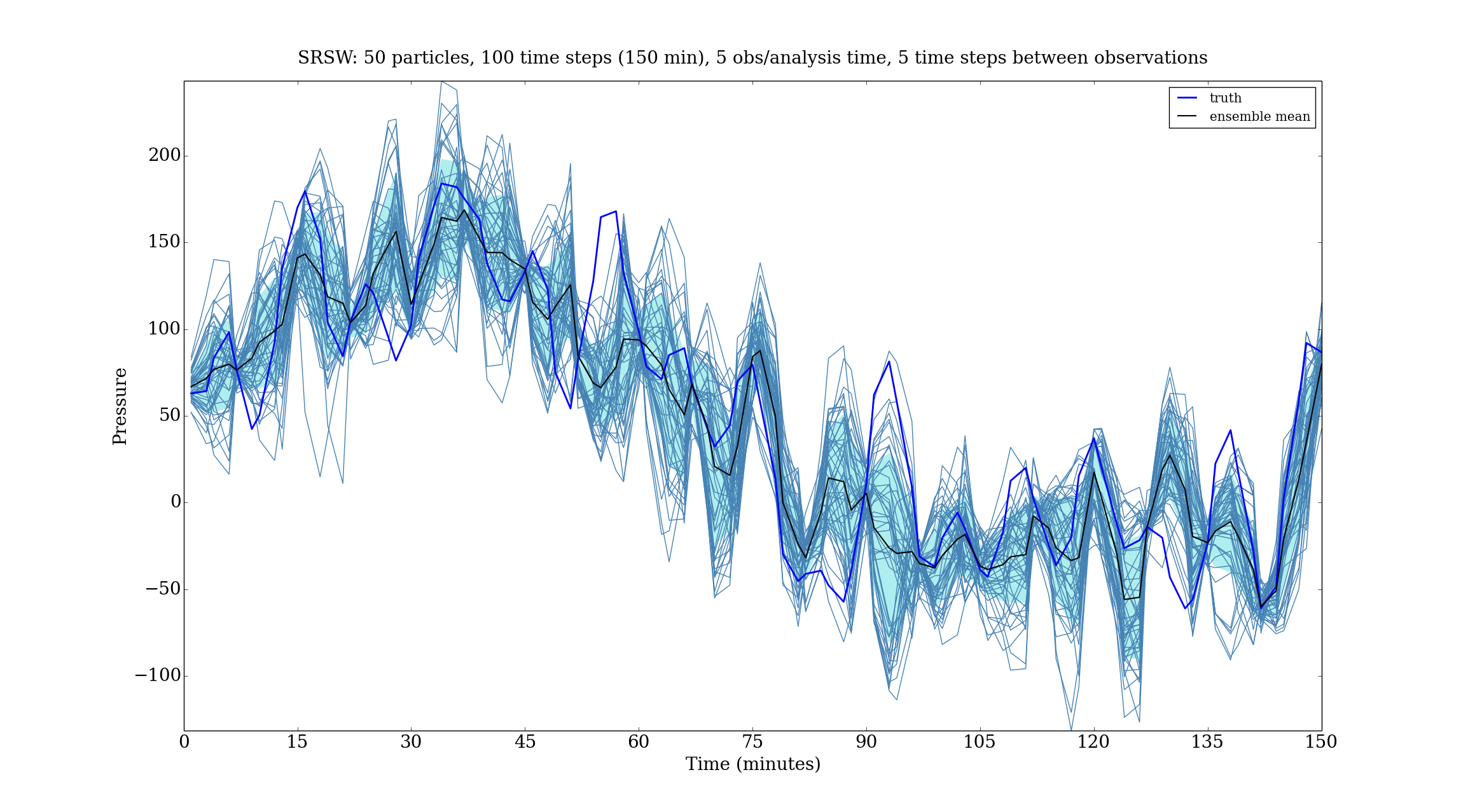}
    \caption{pressure field}
  \label{dasrswS32p}
  \end{subfigure}%
 \begin{subfigure}{.5\linewidth}
    \centering
   \includegraphics[width = \linewidth]{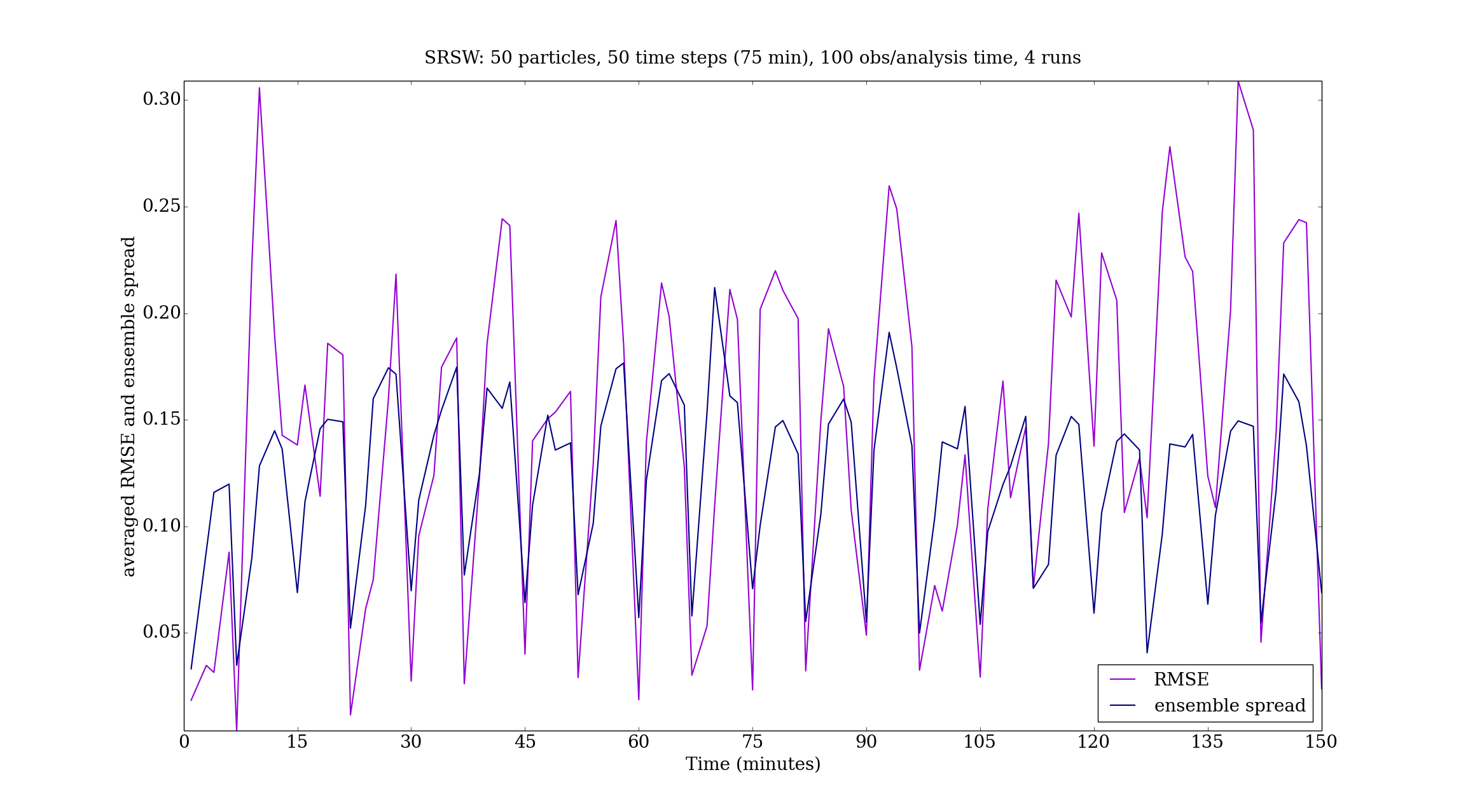}
   \caption{RMSE}
   \label{dasrswS32rmse}
  \end{subfigure}%
  \caption{Evolution of the SRSW model for 100 time steps. The system is observed every 5 time steps, with 5 observations/analysis time. { The cloud of particles follows the truth most of the time.}}
\end{figure}

 We show below { (Figures \ref{dasrswS32p}-\ref{dasrswS32rmse})} the evolution of the SRSW model for 100 time steps, when observed every 5 time steps, with 5 observations at each analysis time. Although the particles have a constant tendency to diverge from the truth, especially after the first assimilation step, by observing the system quite often and making the observations informative enough, we capture the truth most of the time. 
 
\section{Conclusions and further work}

We conclude that our methodology based on a bootstrap particle filter with tempering and jittering can be successfully used to assimilate data for the Lorenz '63 model. The performance of the data assimilation methodology is shown to be influenced by the following factors: 
\begin{itemize}
\item Observational uncertainty: we test this for three values of the observation error (0.1, 0.5, and 5); the results improve as the observation error is reduced. 
\item Initial uncertainty: we test this in two cases (initially this parameter is set to 1, and then increased to 3). The results improve as the initial uncertainty is reduced, although that advantage decreases over time as the initial uncertainty is forgotten. 
\item Linearity of the observation operator - three observation operators are used: a fully linear one ($\mathscr{H}(x,y,z)=(x,y,z)$), a partially linear one ($\mathscr{H}(x,y,z)=(x^2,y,z)$), and a fully nonlinear one ($\mathscr{H}(x,y,z)=(x^2,y^2,z^2)$). We obtain the best results in the linear case. However, it is worth highlighting that we obtain good results also in the nonlinear case, where we see a bimodal distribution. In contrast to this, a standard Kalman filter would probably fail, due to the difficulties generated by multiple modes. 
\end{itemize}
In future research 
we intend to explore also the influence of other parameters such as: the number of particles, the size of data assimilation window, the $ess$ threshold, the model error. 

We conclude that our particle filter based on tempering and jittering can be successfully used also for the SRSW model. We highlight that this particle filter has been designed without using any other approximation methods (such as localisation). According to our results, there are a couple of parameters which have a relevant influence on the proficiency of the particle filter:
\begin{itemize}
\item DA time window: the results are better when we observe the system every 5 time steps, compared to when we observe it every 10 time steps.   
\item Number of observations per assimilation time: we test this for 1 observation, 5 observations, and 100 observations, respectively, per analysis time. The best output is obtained in the last case. However, by reducing the number of observations from 100 to 5, but observing the system every 5 time steps, the truth can be well tracked for a long period of time (Figure \ref{dasrswS32p}).  
\item Model error: the particle filter is shown to be robust with respect to increases in model error. 
\end{itemize}
The experiments show that the system is underobserved. This is to be expected, as the number of degrees of freedom is very high ($\sim 133,440$). Nonetheless, the findings in this paper are sufficiently promising to encourage further in-depth investigations based on this particle filter. 
We intend to explore also the influence of other parameters such as: the number of particles, observational and initial uncertainty, the $ess$ threshold, variations of the stochastic forcing. Moreover, we intend to use this work as a stepping stone for modelling a slice of the real atmosphere using pressure observations collected by DWD using commercial aircraft.

\end{document}